\documentclass[number,seceqn,secthm,dvips]{arxbj}
\usepackage{upgreek}
\usepackage{graphicx}

\aid{0}
\volume{16}
\issue{3}
\pubyear{2010}
\firstpage{858}
\lastpage{881}
\doi{10.3150/09-BEJ235}

\makeatletter
\newcommand{\eqref}[1]{(\ref{#1})}
\newremark{rem}{Remark}[section]
\makeatother

\begin{document}
\begin{frontmatter}

\title{Fractional pure birth processes}
\runtitle{Fractional birth processes}

\begin{aug}
\author[a]{\fnms{Enzo} \snm{Orsingher}\corref{}\thanksref{e1}\ead[label=e1,mark]{enzo.orsingher@uniroma1.it}}
\and
\author[a]{\fnms{Federico} \snm{Polito}\thanksref{e2}\ead[label=e2,mark]{federico.polito@uniroma1.it}}
\runauthor{E. Orsingher and F. Polito}
\address[a]{Dipartimento di Statistica, Probabilit\`a e Stat. Appl.,
``Sapienza'' Universit\`a di Roma, pl.~A.~Moro~5,\\ 00185~Rome, Italy.
E-mails: \printead*{e1,e2}}
\end{aug}

% HISTORY:
\received{\smonth{11} \syear{2008}}
\revised{\smonth{8} \syear{2009}}

% ABSTRACT
%
\begin{abstract}
We consider a fractional version of the classical nonlinear birth
process of which the Yule--Furry model is a particular
case. Fractionality is obtained by replacing the first order time
derivative in the difference-differential
equations which govern the probability law of the process with the
Dzherbashyan--Caputo fractional derivative.
We derive the probability distribution of the number $ \mathcal{N}_\nu
( t )$ of individuals at an
arbitrary time $t$. We also present an interesting representation for
the number of individuals at time $t$,
in the form of the subordination relation $\mathcal{N}_\nu( t
) = \mathcal{N}
( T_{2 \nu} ( t ) )$, where $\mathcal{N} (
t )$ is the classical generalized birth process and $T_{2 \nu}
( t )$ is
a random time whose distribution is related to the fractional
diffusion equation. The fractional linear birth
process is examined in detail in Section 3 and various
forms of its distribution are given and
discussed.
\end{abstract}

% KEYWORDS
%

\end{frontmatter}
%

%s1 ###
\section{Introduction}
We consider a birth process and denote by $\mathcal{N} ( t
)$, $t>0$, the number of components in a
stochastically developing population at time $t$.
Possible examples are the number of particles produced in a
radioactive disintegration and the number
of particles in a cosmic ray shower where death is not permitted.
The probabilities $p_k ( t ) = \operatorname{Pr} \{ \mathcal
{N} ( t ) = k \}$ satisfy the
difference-differential equations
%
%e1.1 ###
\begin{equation}
\label{basic-equation}
\frac{\mathrm{d} p_k}{\mathrm{d}t} = - \lambda_k p_k + \lambda_{k-1} p_{k-1},\qquad   k
\geq1,
\end{equation}
where, at time $t=0$,
%
%e1.2 ###
\begin{equation}
p_k ( 0 ) = \cases{
1,  & \quad  $k=1,$ \cr
0, & \quad   $k \geq2.$}
\end{equation}
This means that we initially have one progenitor igniting the branching
process. For information on this
process, consult Gikhman and Skorokhod \cite{Skoro}, page 322.

Here, we will examine a fractional version of the birth process where
the probabilities are governed by
%
%e1.3 ###
\begin{equation}
\label{fra}
\frac{\mathrm{d}^\nu p_k}{\mathrm{d}t^\nu} = - \lambda_k p_k + \lambda_{k-1} p_{k-1},
  \qquad k \geq1,
\end{equation}
and where the fractional derivative is understood in the
Dzherbashyan--Caputo sense, that is, as
%
%e1.4 ###
\begin{equation}
\label{caputo}
\frac{\mathrm{d}^\nu p_k}{\mathrm{d}t^\nu} = \frac{1}{\Gamma( 1- \nu)} \int
_0^t \frac{({\mathrm{d}}/{\mathrm{d}s}) p_k (
s )}{( t-s )^\nu}\, \mathrm{d}s\qquad  \mbox{for $ 0 < \nu< 1 $}
\end{equation}
(see Podlubny \cite{podlubny}). The use of a Dzherbashyan--Caputo derivative is
preferred because in this case, initial conditions can be expressed in
terms of integer-order derivatives.

Extensions of continuous-time point processes like the homogeneous
Poisson process to the
fractional case have been considered in Jumarie \cite{jumarie}, Cahoy \cite
{cahoy}, Laskin \cite{laskin}, Wang and Wen \cite{xiao1},
Wang, Wen and Zhang \cite{xiao2}, Wang, Zhang and Fan \cite{wang},
Uchaikin and Sibatov \cite{sibatov}, Repin and Saichev \cite{repin}
and Beghin and Orsingher \cite
{orsbeg}. A recently
published paper (Uchaikin, Cahoy and Sibatov \cite{sib})
considers a fractional version of the
Yule--Furry process
where the mean value $\mathbb{E} N_\nu( t )$ is analyzed.

By recursively solving equation (\ref{fra}) (we write $p_k (
t )$, $t>0$, in equations (\ref{fra}) and $p_k^\nu
( t )$ for the solutions), we obtain that
%
%e1.5 ###
\begin{eqnarray}
\label{recurs}
p_k^\nu( t )& =& \operatorname{Pr} \{ \mathcal{N}_\nu( t
) = k \}
\nonumber
\\[-8pt]
\\[-8pt]
\nonumber
& =& \cases{
 \displaystyle\prod _{j=1}^{k-1} \lambda_j  \sum_{m=1}^k
\biggl\{
\frac{1}{\prod_{ l=1, l \neq m }^k
( \lambda_l - \lambda_m )} E_{\nu, 1} ( - \lambda_m
t^\nu) \biggr\}, & \quad $k > 1,$ \vspace*{2pt}\cr
E_{\nu,1}  ( - \lambda_1 t^\nu), &\quad  $k=1.$}
\end{eqnarray}
Result (\ref{recurs}) generalizes the classical distribution of the
birth process (see
Gikhman and Skorokhod \cite{Skoro}, page 322, or Bartlett \cite{bartlett}, page 59), where, instead
of the
exponentials, we have the Mittag--Leffler functions, defined as
%
%e1.6 ###
\begin{equation}
\label{mittag}
E_{\nu,1} ( x ) = \sum_{h=0}^\infty\frac{x^h}{\Gamma
( \nu h+1 )},\qquad
x \in\mathbb{R},  \nu> 0 .
\end{equation}

% The fractional pure bith process has some specific features entailed
%by the fractional derivatives
% \eqref{caputo}. This implies that the evolution of the probabilitites
% $p_k^\nu( t ) = \operatorname{Pr} \{ N_\nu( t )
%= k \}$
% display a global memory and the structure of the inter-times between
%births is accelerated
% as the order of fractionality decreases. We show below that $
% t )$ and $\mathbb{V} $ar $N_\nu( t )$ are
%increasing as the order
% of fractionality $\nu$ decreases. This makes the model presented
%capable of
% representing highly accelerated epidemics and cosmic showers. This is
%similar to
% what has been pointed out for fractional Poisson process (the
%inter-times
% between events have a Mittag-Leffler structure) which is used in
%optics,
% finance and economics as \cite{cahoy} indicates at page 20.

The fractional pure birth process has some specific features entailed
by the
fractional derivative appearing in \eqref{caputo}, which is a
non-local operator.
The process governed by fractional equations (and therefore the related
probabilities $p_k^\nu( t ) = \operatorname{Pr} \{ N_\nu
( t )
= k \},   k \geq1$) displays a slowly decreasing memory which
seems a
characteristic feature of all real systems (for example, the hereditariety
and the related aspects observed in phenomena such as metal fatigue,
magnetic hysteresis and others). Fractional equations of various types
have proven to be
useful in representing different phenomena in optics (light propagation
through random media), transport of charge carriers and also in economics
(a survey of applications can be found in Podlubny \cite{podlubny}).
Below, we show that for the linear birth process $N_\nu( t
),   t>0,$
the mean values $\mathbb{E} N_\nu( t )$,
$\operatorname{\mathbb{V}ar}N_\nu( t )$ are increasing functions
as the order of fractionality $\nu$ decreases. This shows that the
fractional birth process is capable of representing explosively
developing epidemics, accelerated cosmic showers and, in general,
very rapidly expanding populations. This is a feature which
the fractional pure birth process shares with its Poisson
fractional counterpart whose practical applications have been
studied in recent works (see, for example, Laskin \cite{laskin}
and Cahoy \cite{cahoy}).

We are able to show that the fractional birth process $\mathcal{N}_\nu
( t )$ can be represented
as
%
%e1.7 ###
\begin{equation}
\label{iterated-repr}
\mathcal{N}_\nu( t ) = \mathcal{N} ( T_{2 \nu} (
t ) ), \qquad  t>0,   0<\nu\leq1,
\end{equation}
where $T_{2 \nu} ( t )$, $t>0$, is the random time process
whose distribution at time $t$ is obtained
from the fundamental solution to the fractional diffusion equation
(the fractional derivative is defined in~\eqref{caputo})
%
%e1.8 ###
\begin{equation}
\label{funda}
\dfrac{\partial^{2 \nu} u}{\partial t^{2 \nu}} = \dfrac{\partial^2
u}{\partial s^2},\qquad
 0 < \nu\leq1,
\end{equation}
subject to the initial conditions $u ( s,0 ) = \delta(
s )$
for $0<\nu\leq1$ and also $u_t ( s,0 ) = 0$ for $1/2 <
\nu\leq1$,
as
%
%e1.9 ###
\begin{equation}
\operatorname{Pr} \{ T_{2 \nu} ( t ) \in \mathrm{d}s \} = \cases{
2 u_{2 \nu} ( s , t ) \,\mathrm{d}s & \quad $\mbox{for } s>0,$ \cr
0 & \quad $\mbox{for }s<0.$}
\end{equation}
This means that the fractional birth process is a classical birth
process with a random time $T_{2 \nu}
( t )$ which is the sole component of (\ref{iterated-repr})
affected by the fractional derivative. In equation \eqref{funda} and
throughout the whole paper, the fractional derivative must be
understood in
the Dzherbashyan--Caputo sense \eqref{fra}.
The representation (\ref{iterated-repr}) leads to
%
%e1.10 ###
\begin{equation}
\label{other-repr}
\operatorname{Pr} \{ \mathcal{N}_\nu( t ) = k \} = \int
_0^\infty\operatorname{Pr} \{
\mathcal{N} ( s ) = k \} \operatorname{Pr} \{ T_{2 \nu}
( t ) \in \mathrm{d}s \},
\end{equation}
where
%
%e1.11 ###
\begin{equation}
\label{non-linear-classic}
\operatorname{Pr} \{ \mathcal{N} ( s ) = k \} =
\cases{
 \displaystyle\prod _{j=1}^{k-1} \lambda_j
 \sum_{m=1}^k \frac{ \mathrm{e}^{- \lambda_m s}}{ \prod
 _{l=1, l \neq m}^k ( \lambda_l - \lambda_m ) }, & \quad $k
> 1$, $s>0,$\vspace*{2pt} \cr
\mathrm{e}^{-\lambda_1 s}, &\quad $k=1$, $s>0.$}
\end{equation}
Formula (\ref{other-repr}) immediately shows that $\sum_k \operatorname{Pr}
\{ \mathcal{N}_{\nu} ( t )
= k \} = 1$
if and only if\break $\sum_k \operatorname{Pr} \{ \mathcal{N} ( t )
= k \} =1$.
It is well known that the process $\mathcal{N} ( t )$, $ t
> 0$, is such that $\operatorname{Pr}
( \mathcal{N} ( t )
< \infty) = 1 $ for all $t>0$ (non-exploding) if $\sum_k \lambda
_k^{-1} = \infty$ (see
Feller \cite{feller1}, page 452).

A special case of the above fractional birth process is the fractional
linear birth process
where $\lambda_k = \lambda k$. In this case, the distribution (\ref
{recurs}) reduces to the simple
form
%
%e1.12 ###
\begin{equation}
\label{dist}
p_k^\nu( t ) = \sum_{j=1}^{k} \pmatrix{{k-1}\cr{j-1}} ( -1
)^{j-1} E_{\nu,1}
( - \lambda j t^\nu), \qquad k \geq1, t>0.
\end{equation}
For $\nu=1$, we retrieve from (\ref{dist}) the classical geometric
structure of the linear birth process with a
single progenitor, that is,
%
%e1.13 ###
\begin{equation}
p_k^1 ( t ) = ( 1-\mathrm{e}^{- \lambda t} )^{k-1} \mathrm{e}^{-
\lambda t}, \qquad k \geq1, t>0 .
\end{equation}
An interesting qualitative feature of the fractional linear birth
process can be extracted from~(\ref{dist}); it permits us
to highlight the dependence of the branching speed on the order of
fractionality $\nu$. We show in Section
\ref{sec-exp} that
%
%e1.14 ###
\begin{equation}
\label{dependence}
\operatorname{Pr} \{ N_\nu( \mathrm{d}t ) = n_0 +1  |
N_\nu( 0 ) = n_0 \}
\sim\frac{\lambda n_0 ( \mathrm{d}t )^\nu}{\Gamma( \nu+ 1
)}
\end{equation}
and this proves that a decrease in the order of fractionality $\nu$
speeds up the reproduction of individuals. We are not
able to generalize (\ref{dependence}) to the case
%
%e1.15 ###
\begin{equation}
\operatorname{Pr} \{ N_\nu( t + \mathrm{d}t ) = n_0 +1  |
N_\nu( t ) = n_0 \}
\end{equation}
because the process we are investigating is not time-homogeneous. For
the fractional linear birth process,
the representation (\ref{iterated-repr}) reduces to the form
%
%e1.16 ###
\begin{equation}
\label{iterated-repr-lin}
N_\nu( t ) = N ( T_{2 \nu} ( t ) ),\qquad
  t>0,   0<\nu\leq1,
\end{equation}
and has an interesting special structure when $\nu= 1 / 2^n$. For
example, for $n=2$, the random time appearing
in (\ref{iterated-repr-lin}) becomes a folded iterated Brownian
motion. This means that
%
%e1.17 ###
\begin{equation}
\label{iterated-bm}
N_{{1}/{4}} ( t ) = N ( \vert\mathcal{B}_1
( \vert\mathcal{B}_2 ( t )
\vert) \vert) .
\end{equation}
Clearly, $| \mathcal{B}_2 ( t ) |$ is a reflecting Brownian
motion starting from zero
and $|\mathcal{B}_1 ( |\mathcal{B}_2 ( t ) | ) |$
is a reflecting iterated Brownian motion. This permits us to write the
distribution of (\ref{iterated-bm}) in the following form:
%
%e1.18 ###
\begin{eqnarray}
&&\operatorname{Pr} \{ N_{{1}/{4}} ( t ) =k  |
N_{{1}/{4}} ( 0 ) = 1 \}
\nonumber
\\[-8pt]
\\[-8pt]
\nonumber
&&\quad= \int_0^\infty( 1-\mathrm{e}^{- \lambda s} )^{k-1} \mathrm{e}^{- \lambda s}
\biggl\{
2^2 \int_0^\infty\frac{\mathrm{e}^{- {s^2}/{(4 \omega)}}}{\sqrt{2 \uppi2
\omega}} \frac{\mathrm{e}^{- {\omega^2}/{
(4t)}}}{\sqrt{2 \uppi2t}} \,\mathrm{d} \omega\biggr\}\, \mathrm{d}s .
\end{eqnarray}
The case $\nu= 1/2^n$ involves the $(n -1 )$-times iterated
Brownian motion
%
%e1.19 ###
\begin{equation}
\mathcal{I}_{n-1} ( t ) = \mathcal{B}_1 ( \vert
\mathcal{B}_2 ( \cdots\vert\mathcal{B}_n
( t ) \vert\cdots)  | )
\end{equation}
with distribution
\begin{eqnarray}
\label{iterated-distr}
&& \operatorname{Pr} \{ \vert\mathcal{B}_1 ( \vert\mathcal{B}_2
( \cdots\vert\mathcal{B}_n
( t ) \vert\cdots)  ) \vert
\in \mathrm{d}s \}
\nonumber
\\[-8pt]
\\[-8pt]
\nonumber
&&\quad =\mathrm{d}s 2^n \int_0^\infty\frac{\mathrm{e}^{- {s^2}/{(4 \omega_1)}}}{\sqrt{4 \uppi
\omega_1}} \,\mathrm{d} \omega_1
\int_0^\infty\frac{\mathrm{e}^{- {\omega_1^2}/{(4 \omega_2)}}}{\sqrt{4 \uppi
\omega_2}}\, \mathrm{d} \omega_2 \cdots
\int_0^\infty\frac{\mathrm{e}^{- {\omega_{n-1}^2}/{(4t)}}}{\sqrt{4 \uppi t}}\, \mathrm{d}
\omega_{n-1} . \nonumber
\end{eqnarray}
For details on this point, see Orsingher and Beghin \cite{ors2008}.

%s2 ###
\section{The distribution function for the generalized fractional
birth process}
We now present the explicit distribution
%
%e2.1 ###
\begin{equation}
\operatorname{Pr} \{ \mathcal{N}_\nu( t )=k  |
\mathcal{N}_\nu ( 0 ) = 1
\} = p_k^\nu( t ),\qquad   t>0,   k \geq1,   0 <
\nu\leq1,
\end{equation}
of the number of individuals in the population expanding according to
(\ref{fra}). Our technique is based on successive
applications of the Laplace transform. Our first result is the
following theorem.
\begin{thm}
\label{theorem-general}
The solution to the fractional equations
%
%e2.2 ###
\begin{equation}
\label{cauchy-again}
\cases{
\displaystyle\frac{\mathrm{d}^\nu p_k}{\mathrm{d} t^\nu} = - \lambda_k p_k + \lambda_{k-1} p_{k-1},
& \quad $k \geq1, 0 < \nu\leq1,$ \cr
p_k ( 0 ) = \cases{
1, &\quad  $k = 1,$ \cr
0, & \quad $k \geq2,$}&}
\end{equation}
is given by
%
%e2.3 ###
\begin{eqnarray}
\label{recurs-again}
p_k^\nu( t ) &=& \operatorname{Pr} \{ \mathcal{N}_\nu( t
) = k |\mathcal{N}_\nu ( 0 ) = 1\}
\nonumber
\\[-8pt]
\\[-8pt]
\nonumber
&=&
\cases{
 \displaystyle\prod _{j=1}^{k-1} \lambda_j  \sum_{m=1}^k
\biggl\{
\frac{1}{\prod _{ l=1, l \neq m }^k
( \lambda_l - \lambda_m )} E_{\nu, 1} ( - \lambda_m
t^\nu) \biggr\}, &\quad $ k>1,$\cr
E_{\nu,1} ( - \lambda_1 t^\nu), &\quad $ k=1$.}
\end{eqnarray}
\end{thm}
\begin{pf}
We prove the result (\ref{recurs-again}) by a recursive procedure.

For $k=1$, the equation
%
%e2.4 ###
\begin{equation}
\frac{\mathrm{d}^\nu p_1}{\mathrm{d}t^\nu} = -\lambda_1 p_1, \qquad   p_1 ( 0 ) =1,
\end{equation}
is immediately solved by
%
%e2.5 ###
\begin{equation}
p_1^\nu( t ) = E_{\nu,1} ( - \lambda_1 t^\nu).
\end{equation}
For $k=2$, equation (\ref{fra}) becomes
%
%e2.6 ###
\begin{equation}
\label{kappa2}
\cases{
\dfrac{\mathrm{d}^\nu p_2}{\mathrm{d}t^\nu} = - \lambda_2 p_2 + \lambda_1 E_{\nu,1}
( -\lambda_1 t^\nu
), \vspace*{2pt}\cr
p_2 ( 0 ) = 0.}
\end{equation}
In view of the fact that
%
%e2.7 ###
\begin{equation}
\label{mittag-laplace}
\int_0^\infty \mathrm{e}^{- \mu t} E_{\nu, 1} ( - \lambda_1 t^\nu)\,
\mathrm{d}t =
\frac{\mu^{\nu-1}}{\mu^\nu+ \lambda_1},
\end{equation}
the Laplace transform of (\ref{kappa2}) yields
%
%e2.8 ###
\begin{equation}
\label{laplace-kappa2}
L_2 ( \mu) =
\frac{\lambda_1 \mu^{\nu-1}}{\lambda_2 - \lambda_1}\biggl [ \frac
{1}{\mu^\nu+
\lambda_1} - \frac{1}{\mu^\nu+ \lambda_2} \biggr] .
\end{equation}
In the light of (\ref{mittag-laplace}), from (\ref{laplace-kappa2}), we
can determine the probability
$p_2^\nu( t )$:
%
%e2.9 ###
\begin{equation}
p_2^\nu( t ) = [ E_{\nu,1} ( - \lambda_1 t^\nu
) -
E_{\nu,1} ( - \lambda_2 t^\nu) ] \frac{\lambda_1}{
\lambda_2 - \lambda_1}.
\end{equation}
Now, the Laplace transform of
%
%e2.10 ###
\begin{equation}
\frac{\mathrm{d}^\nu p_3}{\mathrm{d}t^\nu} = - \lambda_3 p_3 + \frac{\lambda_2 \lambda
_1}{\lambda_2
- \lambda_1} [ E_{\nu,1} ( - \lambda_1 t^\nu) - E_{\nu
, 1}
( -\lambda_2 t^\nu) ]
\end{equation}
yields, after some computation,\vspace*{2pt}
%
%e2.11 ###
\begin{eqnarray}
&&L_3 ( \mu) = \lambda_2 \lambda_1 \mu^{\nu-1} \biggl[ \frac{1}{( \lambda_2
- \lambda_1
) ( \lambda_3 - \lambda_1 )} \frac{1}{\mu^\nu+
\lambda_1}\nonumber\\[2pt]
&&\phantom{L_3 ( \mu) = \lambda_2 \lambda_1 \mu^{\nu-1} \biggl[}{} + \frac{1}{( \lambda_1 -\lambda_2 ) ( \lambda
_3 - \lambda_2 )}
\frac{1}{\mu^\nu+ \lambda_2} \\[2pt]
&&\phantom{L_3 ( \mu) = \lambda_2 \lambda_1 \mu^{\nu-1} \biggl[}{}+ \frac{1}{( \lambda_1 - \lambda_3
) (
\lambda_2 - \lambda_3 )} \frac{1}{\mu^\nu+ \lambda_3} \biggr] .\vspace*{2pt}
\nonumber
\end{eqnarray}
From this result, it is clear that\vspace*{1pt}
%
%e2.12 ###
\begin{eqnarray}
&&p_3^\nu( t ) = \lambda_2 \lambda_1 \biggl[ \frac
{1}{( \lambda_2 - \lambda_1
) ( \lambda_3 - \lambda_1 )} E_{\nu,1} ( -
\lambda_1 t^\nu)\nonumber\\[2pt]
 &&\phantom{p_3^\nu( t ) = \lambda_2 \lambda_1 \biggl[}{}+ \frac{1}{( \lambda_1 -\lambda_2 ) ( \lambda
_3 - \lambda_2 )}
E_{\nu,1} ( - \lambda_2 t^\nu) \\[2pt]
 &&\phantom{p_3^\nu( t ) = \lambda_2 \lambda_1 \biggl[}{}+ \frac{1}{( \lambda_1
- \lambda_3 ) (
\lambda_2 - \lambda_3 )} E_{\nu,1} ( - \lambda_3 t^\nu
) \biggr] . \nonumber\vspace*{2pt}
\end{eqnarray}
The procedure for $k > 3$ becomes more complicated. However, the
special case $k=4$ is instructive and
so we treat it first.

The Laplace transform of the equation\vspace*{1pt}
%
%e2.13 ###
\begin{eqnarray}
&&\frac{\mathrm{d}^\nu p_4}{\mathrm{d}t^\nu} =  - \lambda_4 p_4 + \lambda_1 \lambda_2
\lambda_3
\biggl[ \frac{1}{( \lambda_2 - \lambda_1 ) ( \lambda_3
- \lambda_1 )}
E_{\nu,1} ( - \lambda_1 t^\nu) \nonumber\\[2pt]
&&\phantom{\frac{\mathrm{d}^\nu p_4}{\mathrm{d}t^\nu} =  - \lambda_4 p_4 + \lambda_1 \lambda_2
\lambda_3
\biggl[}{} +  \frac{1}{ ( \lambda_1 - \lambda_2 ) (
\lambda_3 - \lambda_2 )}
E_{\nu,1} ( - \lambda_2 t^\nu) \\[2pt]
&&\phantom{\frac{\mathrm{d}^\nu p_4}{\mathrm{d}t^\nu} =  - \lambda_4 p_4 + \lambda_1 \lambda_2
\lambda_3
\biggl[}{}+ \frac{1}{ ( \lambda
_1 - \lambda_3
) ( \lambda_2 - \lambda_3 ) } E_{\nu,1} ( -
\lambda_2 t^\nu
) \biggr], \nonumber\vspace*{2pt}
\end{eqnarray}
subject to the initial condition $p_4 ( 0 ) =0$, becomes\vspace*{2pt}
%
%e2.14 ###
\begin{eqnarray}
\label{laplace-for-4}
&&L_4 ( \mu)
= \lambda_1 \lambda_2 \lambda_3 \mu^{\nu-1} \biggl[ \frac{1}{
(
\lambda_2 - \lambda_1 ) ( \lambda_3 - \lambda_1 )
(
\lambda_4 - \lambda_1 ) } \biggl\{ \frac{1}{\mu^\nu+ \lambda_1} -
\frac{1}{\mu^\nu+ \lambda_4} \biggr\} \nonumber\\[2pt]
&&\phantom{L_4 ( \mu)
= \lambda_1 \lambda_2 \lambda_3 \mu^{\nu-1} \biggl[}{}+ \frac{1}{( \lambda_1 - \lambda_2 ) ( \lambda_3 -
\lambda_2 )
( \lambda_4 - \lambda_2 ) } \biggl\{ \frac{1}{\mu^\nu+
\lambda_2} -
\frac{1}{\mu^\nu+ \lambda_4} \biggr\} \\[2pt]
&&\phantom{L_4 ( \mu)
= \lambda_1 \lambda_2 \lambda_3 \mu^{\nu-1} \biggl[}{}  + \frac{1}{( \lambda_1 - \lambda_3 ) ( \lambda
_2 - \lambda_3 )
( \lambda_4 - \lambda_3 )} \biggl\{ \frac{1}{\mu^\nu+
\lambda_3 } -
\frac{1}{\mu^\nu+ \lambda_4 } \biggr\} \biggr] .\qquad  \nonumber
\end{eqnarray}
The critical point of the proof is to show that
%
%e2.15 ###
\begin{eqnarray}
&&- \bigl( ( \lambda_3 - \lambda_2 ) ( \lambda_4 -
\lambda_2 )
( \lambda_4 - \lambda_3 ) - ( \lambda_3 - \lambda_1
)
( \lambda_4 - \lambda_1 ) ( \lambda_4 - \lambda_3
) \nonumber\\
&&\qquad{}+
( \lambda_2 - \lambda_1 ) ( \lambda_4 - \lambda_1
)
( \lambda_4 - \lambda_2 )\bigr)\nonumber
\\[-8pt]
\\[-8pt]
\nonumber
&&\qquad {}\times\frac1{( \lambda_2 - \lambda_1
)
( \lambda_3 - \lambda_1 ) ( \lambda_4 - \lambda_1
)
( \lambda_3 - \lambda_2 ) ( \lambda_4 - \lambda_2
)
( \lambda_4 - \lambda_3 )}  \\
&&\quad = \frac{1}{ ( \lambda_1 - \lambda_4 )
( \lambda_2 - \lambda_4 ) ( \lambda_3 - \lambda_4
)} . \nonumber
\end{eqnarray}
We note that
%
%e2.16 ###
\begin{eqnarray}
\label{vandermonde}
0 &=&  \det\pmatrix{
1 & 1 & 1 & 1 \cr
1 & 1 & 1 & 1 \cr
\lambda_1 & \lambda_2 & \lambda_3 & \lambda_4 \vspace*{2pt}\cr
\lambda_1^2 & \lambda_2^2 & \lambda_3^2 & \lambda_4^2}\nonumber\\
& =& \det\pmatrix{
1 & 1 & 1 \cr
\lambda_2 & \lambda_3 & \lambda_4\vspace*{2pt} \cr
\lambda_2^2 & \lambda_3^2 & \lambda_4^2}\nonumber\\
&&{} - \det\pmatrix{
1 & 1 & 1 \cr
\lambda_1 & \lambda_3 & \lambda_4 \vspace*{2pt}\cr
\lambda_1^2 & \lambda_3^2 & \lambda_4^2} + \det\pmatrix{
1 & 1 & 1 \cr
\lambda_1 & \lambda_2 & \lambda_4 \vspace*{2pt}\cr
\lambda_1^2 & \lambda_2^2 & \lambda_4^2} - \det\pmatrix{
1 & 1 & 1 \cr
\lambda_1 & \lambda_2 & \lambda_3 \vspace*{2pt}\cr
\lambda_1^2 & \lambda_2^2 & \lambda_3^2}\\
& =& ( \lambda_3 - \lambda_2 ) ( \lambda_4 - \lambda
_2 )
( \lambda_4 - \lambda_3 ) - ( \lambda_3 - \lambda_1
)
( \lambda_4 - \lambda_1 ) ( \lambda_4 - \lambda_3
) \nonumber\\
&&{} + ( \lambda_2 - \lambda_1 ) ( \lambda_4 - \lambda_1
)
( \lambda_4 - \lambda_2 )
- ( \lambda_2 - \lambda_1 ) ( \lambda_3 - \lambda_1
)
( \lambda_3 - \lambda_2 ), \nonumber
\end{eqnarray}
where, in the last step, the Vandermonde formula is applied.

By inserting (\ref{vandermonde}) into (\ref{laplace-for-4}), we now
have that
%
%e2.17 ###
\begin{eqnarray}
\label{kappa4}
&&L_4 ( \mu) =  \lambda_1 \lambda_2 \lambda_3 \mu^{\nu
-1} \biggl[
\frac{1}{( \lambda_2 - \lambda_1 ) ( \lambda_3 -
\lambda_1 )
( \lambda_4 - \lambda_1 )} \frac{1}{ \mu^\nu+ \lambda_1} \nonumber\\
&&\phantom{L_4 ( \mu) =  \lambda_1 \lambda_2 \lambda_3 \mu^{\nu
-1} \biggl[}{}+ \frac{1}{( \lambda_1 - \lambda_2 ) ( \lambda_3 -
\lambda_2 )
( \lambda_4 - \lambda_2 ) } \frac{1}{ \mu^\nu+ \lambda_2 }
\nonumber
\\[-8pt]
\\[-8pt]
\nonumber
&&\phantom{L_4 ( \mu) =  \lambda_1 \lambda_2 \lambda_3 \mu^{\nu
-1} \biggl[}{}+\frac{1}{( \lambda_1 - \lambda_3 ) ( \lambda_2 -
\lambda_3 )
( \lambda_4 - \lambda_3 ) } \frac{1}{\mu^\nu+ \lambda_3 }\\
&&\phantom{L_4 ( \mu) =  \lambda_1 \lambda_2 \lambda_3 \mu^{\nu
-1} \biggl[}{}+ \frac{1}{( \lambda_1 - \lambda_4 ) ( \lambda
_2 -
\lambda_4 ) ( \lambda_3 - \lambda_4 ) } \frac{1}{\mu
^\nu+ \lambda_4}
\biggr] \nonumber
\end{eqnarray}
so that by inverting (\ref{kappa4}), we obtain the following result:
%
%e2.18 ###
\begin{equation}
p_4^\nu( t ) =  \prod _{j=1}^3 \lambda_j
\Biggl\{
 \sum_{m=1}^4 \frac{1}{\prod _{ l=1, l
\neq m }^4
( \lambda_l - \lambda_m )} E_{\nu, 1} ( - \lambda_m
t^\nu) \Biggr\} .
\end{equation}
We now tackle the problem of showing that (\ref{recurs-again}) solves
the Cauchy problem
(\ref{cauchy-again}) for all $k >1$, by induction. This means that we
must solve
%
%e2.19 ###
\begin{equation}
\label{cauchy-for-k}
\cases{
\dfrac{\mathrm{d}^\nu p_k}{\mathrm{d}t^\nu} = - \lambda_kp_k +
 \displaystyle\prod _{j=1}^{k-1} \lambda_j \Biggl\{  \sum
_{m=1}^{k-1}
\frac{1}{\prod _{ l=1, l \neq m }^{k-1}
( \lambda_l - \lambda_m )} E_{\nu, 1} (
- \lambda_m t^\nu) \Biggr\}, \cr
p_k ( 0 ) = 0,}\qquad
  k > 4.
\end{equation}
The Laplace transform of (\ref{cauchy-for-k}) reads
%
%e2.20 ###
\begin{eqnarray}
\label{last-term}
&&L_k ( \mu) =
 \prod _{j=1}^{k-1} \lambda_j \Biggl[
 \sum_{m=1}^{k-1} \frac{\mu^{\nu-1}}{\prod
_{ l=1, l \neq m }^{k}
( \lambda_l - \lambda_m )} \frac{1}{\mu^\nu+ \lambda_m}
\nonumber
\\[-8pt]
\\[-8pt]
\nonumber
&&\phantom{L_k ( \mu) =
 \prod _{j=1}^{k-1} \lambda_j \Biggl[}{} -
\frac{\mu^{\nu-1}}{\mu^\nu+ \lambda_k}  \sum_{m=1}^{k-1}
\frac{1}{\prod _{l=1, l \neq m }^{k}
( \lambda_l - \lambda_m )} \Biggr] .
\end{eqnarray}
We must now prove that
%
%e2.21 ###
\begin{equation}
\label{relation}
-  \sum_{m=1}^{k-1} \frac{1}{\prod _{ l=1, l \neq m }^{k}
( \lambda_l - \lambda_m )} = \frac{1}{\prod
_{ l=1, l \neq k }^{k}
( \lambda_l - \lambda_k )}
\end{equation}
and this relation is also important for the proof of (\ref{non-linear-classic}).

In order to prove (\ref{relation}), we rewrite the left-hand side as
%
%e2.22 ###
\begin{equation}
\label{rewrite}
-  \sum_{m =1}^{k-1} \frac{  \prod _{h=1}^{k-1}
 \prod _{l >h}^k ( \lambda_l - \lambda_h ) }{
 \prod _{l=1, l \neq m}^k ( \lambda_l
- \lambda_m )}
\cdot\frac{1}{  \prod _{h=1}^{k-1}
 \prod _{l>h}^k ( \lambda_l - \lambda_h ) }
\end{equation}
and concentrate our attention on the numerator of (\ref{rewrite}). By
analogy with the calculations
in~(\ref{vandermonde}), we have that
%
%e2.23 ###
\begin{eqnarray}
\label{vander-general}
0 &=& \det\pmatrix{
1 & 1 & \cdots& 1 & \cdots& 1 \cr
1 & 1 & \cdots& 1 & \cdots& 1 \cr
\lambda_1 & \lambda_2 & \cdots& \lambda_m & \cdots& \lambda_k \cr
\cdots& \cdots& \cdots& \cdots& \cdots& \cdots\cr
\lambda_1^{k-2} & \lambda_2^{k-2} & \cdots& \lambda_m^{k-2} & \cdots
& \lambda_k^{k-2} }\nonumber\\
%%%%%%%%%%%%%%%%%%%%%%%%%%%%%%%%%%%%%
&= & \sum_{m=1}^k ( -1 )^{m-1} \det\pmatrix{
1 & \cdots& 1 & 1 & \cdots& 1 \cr
\lambda_1 & \cdots& \lambda_{m-1} & \lambda_{m+1} & \cdots& \lambda
_k\vspace*{2pt} \cr
\lambda_1^{k-2} & \cdots& \lambda_{m-1}^{k-2} & \lambda_{m+1}^{k-2}
& \cdots& \lambda_k^{k-2} }\\
%%%%%%%%%%%%%%%%%%%%%%%%%%%%%%%%%%%%%%%%%%%
&= & \sum_{m=1}^{k} \frac{  \prod _{h=1}^{k-1}
 \prod _{l>h}^{k} ( \lambda_l - \lambda_h ) }{
 \prod _{ l=1,l \neq m }^k ( \lambda
_l - \lambda_m ) } =
\sum_{m=1}^{k-1} \frac{  \prod _{h=1}^{k-1}
 \prod _{l>h}^{k} ( \lambda_l - \lambda_h ) }{
 \prod _{ l=1 , l \neq m }^{k} (
\lambda_l - \lambda_m ) } +
\frac{  \prod _{h=1}^{k-1}
 \prod _{l>h}^{k} ( \lambda_l - \lambda_h ) }{
 \prod _{ l=1, l \neq k }^k ( \lambda
_l - \lambda_k ) }
. \nonumber
\end{eqnarray}
In the third step of (\ref{vander-general}), we applied the Vandermonde
formula and considered
the fact that the~$n$th column is missing. It must also be taken into
account that
%
%e2.24 ###
\begin{eqnarray}
&& \frac{  \prod _{l>1}^k ( \lambda_l - \lambda_1
)}{ (
\lambda_m - \lambda_1 )}   \cdot  \frac{  \prod
 _{l>2}^k (
\lambda_l - \lambda_2 )}{ ( \lambda_m - \lambda_2 )}
  \cdots
\frac{  \prod _{l>m-1}^k ( \lambda_l - \lambda
_{m-1} )}{ (
\lambda_m - \lambda_{m-1} )}   \nonumber\\
&&\qquad{}\times  \frac{  \prod _{l>m}^k (
\lambda_l - \lambda_m )}{  \prod _{l>m}^k (
\lambda_l - \lambda_m )}   \cdot   \prod
_{l>m+1}^k (
\lambda_l - \lambda_{m+1} )   \cdots \prod
_{l>k-1}^k (
\lambda_l - \lambda_{k-1} )\\
&&\quad= \frac{ \prod _{h=1}^{k-1}  \prod
_{l>h}^k (
\lambda_l - \lambda_h )}{ ( -1 )^{m-1}
 \prod _{l=1, l \neq m}^k ( \lambda_l
- \lambda_m )} .
\nonumber
\end{eqnarray}
From (\ref{rewrite}) and (\ref{vander-general}), we have that
%
%e2.25 ###
\begin{eqnarray}
\label{rewrite2}
\hspace*{-20pt}-  \sum_{m=1}^{k-1} \frac{1}{\prod _{ l=1, l \neq m }^{k}
( \lambda_l - \lambda_m )} &=& - \sum_{m=1}^{k-1}
\frac{  \prod _{h=1}^{k-1}  \prod
_{l>h}^{k} ( \lambda_l
- \lambda_h )}{  \prod _{l=1,  l \neq
m}^k (
\lambda_l - \lambda_m )} \cdot\frac{1}{  \prod
_{h=1}^{k-1}
 \prod _{l>h}^k ( \lambda_l - \lambda_h )}
 \nonumber
 \\[-8pt]
 \\[-8pt]
 \nonumber
 &=& \frac{1}{  \prod _{ l=1,l \neq k}^k
( \lambda_l - \lambda_k ) }.
\end{eqnarray}
In view of (\ref{rewrite2}), we can write that
%
%e2.26 ###
\begin{equation}
\label{rewrite-lapl}
L_k ( \mu) =  \displaystyle\prod _{j=1}^{k-1} \lambda_j
 \sum_{m=1}^k \frac{\mu^{\nu-1}}{
 \prod _{ l=1,l \neq m }^k
( \lambda_l - \lambda_m )}
\cdot\frac{1}{\mu^\nu+ \lambda_m}
\end{equation}
because the $k$th term of (\ref{rewrite-lapl}) coincides with the last
term of (\ref{last-term})
and therefore, by inversion of the Laplace transform, we get (\ref
{recurs-again}).
\end{pf}

\begin{rem}
We now prove that for the generalized fractional birth process, the
representation
%
%e2.27 ###
\begin{equation}
\label{iterative}
\mathcal{N}_\nu( t ) = \mathcal{N} ( T_{2 \nu} (
t ) ),\qquad
  t>0,   0< \nu\leq1,
\end{equation}
holds. This means that the process under investigation can be viewed as
a generalized birth process
at a random time $T_{2 \nu} ( t )$, $t>0$, whose
distribution is the folded solution to the
fractional diffusion equation (\ref{funda}).
%
%e2.28 ###
\begin{eqnarray}
\label{lapl-nonlin-g}
&& \int_0^\infty \mathrm{e}^{- \mu t} \mathcal{G}_\nu( u,t ) \,\mathrm{d}t
\nonumber\\
&&\quad\stackrel{\mathrm{by\mbox{ } (\ref{recurs-again})}}{=}\int_0^\infty\Biggl\{ \sum_{k=2}^\infty u^k \prod_{j=1}^{k-1} \lambda
_j \sum_{m=1}^k
\frac{ E_{\nu,1} ( - \lambda_m t^\nu) }{  \prod
 _{j \neq m}^k
( \lambda_j - \lambda_m ) } + u E_{\nu,1} ( - \lambda
_1 t^\nu
) \Biggr\} \mathrm{e}^{- \mu t} \,\mathrm{d}t \nonumber\\
 &&\quad= \sum_{k=2}^\infty u^k \prod_{j=1}^{k-1} \lambda_j \sum_{m=1}^k
\frac{\mu^{\nu-1}}{\mu^\nu+ \lambda_m}
\frac{1}{ \prod _{j \neq m}^k ( \lambda_j - \lambda
_m ) } +
\frac{u \mu^{\nu-1}}{\mu^\nu+ \lambda_1}\\
&&\quad = \int_0^\infty\Biggl\{ \sum_{k=2}^\infty u^k \prod_{j=1}^{k-1}
\lambda_j \sum_{m=1}^k \frac{\mu^{\nu-1}}{
 \prod _{j \neq m}^k ( \lambda_j - \lambda_m
) } \mathrm{e}^{- s ( \mu^\nu+
\lambda_m ) } + u \mathrm{e}^{-s ( \mu^\nu+ \lambda_1 )}
\Biggr\} \,\mathrm{d}s
\nonumber\\
&&\quad= \int_0^\infty\mathcal{G} ( u,s ) \mu^{\nu-1} \mathrm{e}^{- s
\mu^\nu} \,\mathrm{d}s
= \int_0^\infty\mathcal{G} ( u,s ) \int_0^\infty \mathrm{e}^{- \mu
t} f_{T_{2 \nu}} (
s,t ) \,\mathrm{d}t \,  \mathrm{d}s \nonumber\\
&&\quad = \int_0^\infty \mathrm{e}^{- \mu t} \biggl\{ \int_0^\infty\mathcal{G} (
u,s ) f_{T_{2 \nu}}
( s,t ) \,\mathrm{d}s \biggr\}\, \mathrm{d}t , \nonumber
\end{eqnarray}
where
%
%e2.29 ###
\begin{equation}
\int_0^\infty \mathrm{e}^{- \mu t } f_{T_{2 \nu}} ( s,t ) \,\mathrm{d}t = \mu
^{\nu-1} \mathrm{e}^{- s \mu^\nu},\qquad
  s >0 ,
\end{equation}
is the Laplace transform of the folded solution to (\ref{funda}).
From (\ref{lapl-nonlin-g}), we infer that
%
%e2.30 ###
\begin{equation}
\mathcal{G}_\nu( u,t ) = \int_0^\infty\mathcal{G} (
u,s ) f_{T_{2 \nu}}
( s,t )\, \mathrm{d}s
\end{equation}
and from this, the representation (\ref{iterative}) follows.
\end{rem}

\begin{rem}
The relation (\ref{iterative}) permits us to conclude that the
functions (\ref{recurs-again})
are non-negative because
%
%e2.31 ###
\begin{equation}
\operatorname{Pr} \{ \mathcal{N}_\nu( t ) = k \} = \int
_0^\infty\operatorname{Pr} \{
\mathcal{N} ( s ) = k \} \operatorname{Pr} \{ T_{2 \nu}
( t ) \in \mathrm{d}s \},
\end{equation}
and $\operatorname{Pr} \{ \mathcal{N} ( s ) = k \} > 0$
and $\sum_k \operatorname{Pr} \{ \mathcal{N}
( s ) = k \} = 1,$ as shown, for example, in Feller \cite
{feller1}, page 452.
Furthermore, the fractional birth process is non-exploding if and only
if $\sum_k ( 1/ \lambda_k ) = \infty$
for all values of $0 < \nu\leq1$.
\end{rem}

%s3 ###
\section{The fractional linear birth process}
\label{sec-exp}
In this section, we examine in detail a special case of the previous
fractional birth process,
namely the fractional linear birth process which generalizes the
classical Yule--Furry model.
The birth rates in this case have the form
%
%e3.1 ###
\begin{equation}
\label{linear-rates}
\lambda_k = \lambda k,\qquad  \mbox{$\lambda> 0$, $k \geq1$},
\end{equation}
and indicate that new births occur with a probability proportional to
the size of the population.
We denote by $N_\nu( t )$ the number of individuals in the
population expanding
according to the rates~(\ref{linear-rates}) and we have that the probabilities
%
%e3.2 ###
\begin{equation}
\label{pknu}
p_k^\nu( t ) = \operatorname{Pr} \{ N_\nu( t ) =
k \vert N_\nu
( 0 ) = 1 \}, \qquad  k \geq1,
\end{equation}
satisfy the difference-differential equations
%
%e3.3 ###
\begin{equation}
\label{diff-diff-eq}
\cases{
\dfrac{\mathrm{d}^\nu p_k}{\mathrm{d} t^\nu} = - \lambda k p_k + \lambda(k-1 )
p_{k-1}, &\quad $0 < \nu\leq1, k \geq1$, \cr
p_k ( 0 ) = \cases{
1, &\quad  $k = 1,$ \cr
0, &\quad $ k \geq2.$ }&}
\end{equation}
The distribution (\ref{pknu}) can be obtained as a particular case of
(\ref{recurs-again}) or directly, by means
of a completely different approach, as follows.
\begin{thm}
The distribution of the fractional linear birth process with a simple
initial progenitor has the form
%
%e3.4 ###
\begin{eqnarray}
\label{Pk}
p_k^\nu( t ) &=& \operatorname{{\rm Pr}} \{ N_\nu( t
) = k
| N_\nu( 0 ) = 1 \}
\nonumber
\\[-8pt]
\\[-8pt]
\nonumber
&=&
\sum_{j=1}^k \pmatrix{{k-1}\cr{j-1}} ( -1 )^{j-1} E_{\nu,1} (
-\lambda
j t^\nu), \qquad  k \geq1,   0< \nu\leq1,
\end{eqnarray}
where $E_{\nu,1} ( x )$ is the Mittag--Leffler function (\ref
{mittag}).
\end{thm}

\begin{pf}
We can prove the result (\ref{Pk}) by solving equation (\ref
{diff-diff-eq}) recursively. This means that
$p_{k-1}^\nu( t )$ has the form (\ref{Pk}), so $p_k^\nu
( t )$ maintains the same structure. This is tantamount
to solving the Cauchy problem
%
%e3.5 ###
\begin{equation}
\label{problem}
\cases{
\dfrac{\mathrm{d}^\nu p_k ( t ) }{\mathrm{d}t^\nu} = - \lambda k p_k ( t
)
+ \lambda( k-1 )  \displaystyle\sum_{j=1}^{k-1} \pmatrix{{k-2}\cr{j-1}}
( -1 )^{j-1}
E_{\nu,1} (- \lambda j t^\nu),\vspace*{2pt}\cr
p_k ( 0 ) = 0, \qquad  k > 1.}
\end{equation}
By applying the Laplace transform
$
L_{k, \nu} ( \mu) = \int_0^\infty \mathrm{e}^{- \mu t} p_k ( t
) \,\mathrm{d}t
$
to (\ref{problem}), we have that
%
%e3.6 ###
\begin{equation}
\label{laplace}
L_{k,\nu} ( \mu) = \lambda( k-1 )
\Biggl\{ \sum_{j=1}^{k-1} \pmatrix{{k-2}\cr{j-1} }(
-1 )^{j-1} \frac{\mu^{k-1}}{\mu^\nu+ \lambda j} \Biggr\} \frac
{1}{\mu^\nu+ \lambda k} .
\end{equation}
Conveniently, the Laplace transform (\ref{laplace}) can be written as
%
%e3.7 ###
\begin{eqnarray}
\label{expansion}
 L_{k,\nu} ( \mu) &=& \mu^{\nu-1} \biggl\{ \biggl[ \frac
{1}{\mu^\nu+ \lambda}
- \frac{1}{\mu^\nu
+ \lambda k } \biggr] - ( k-1 ) \biggl[ \frac{1}{\mu^\nu+ 2
\lambda} - \frac{1}{\mu^\nu
+ \lambda k} \biggr] \nonumber\\
&&\hspace*{28pt}{} + \frac{ ( k-1 ) ( k-2 ) }{2} \biggl[ \frac
{1}{\mu^\nu+ 3 \lambda} -
\frac{1}{\mu^\nu+ k \lambda} \biggr] + \cdots
\nonumber
\\[-8pt]
\\[-8pt]
\nonumber
&&\hspace*{28pt}{}+ ( k-1 ) ( -1 )^{k-2} \biggl[
\frac{1}{\mu^\nu+ ( k-1 ) \lambda} - \frac{1}{\mu^\nu+
\lambda k}
\biggr] \biggr\}\\
& =& \mu^{\nu-1} \sum_{j=1}^{k-1} \pmatrix{{k-1}\cr{j-1}} ( -1
)^{j-1} \frac{1}{\mu^\nu+
j \lambda} - \frac{\mu^{\nu-1}}{\mu^\nu+ \lambda k} \sum_{j=1}^{k-1}
\pmatrix{{k-1}\cr{j-1}} ( -1
)^{j-1} .\nonumber
\end{eqnarray}
This permits us to conclude that
%
%e3.8 ###
\begin{equation}
\label{laplace-k}
L_{k,\nu} ( \mu) = \mu^{\nu-1} \sum_{j=1}^k \pmatrix{
{k-1}\cr{j-1}} ( -1 )^{j-1}
\frac{1}{\mu^\nu+ j \lambda} .
\end{equation}
By inverting (\ref{laplace-k}), we immediately arrive at the result
(\ref{Pk}).
\end{pf}
For $\nu=1$, (\ref{laplace-k}) can be written as
%
%e3.9 ###
\begin{eqnarray}
\int_0^\infty \mathrm{e}^{- \mu t} p_k^1 ( t ) \,\mathrm{d}t & =& \int
_0^\infty \mathrm{e}^{- \lambda t} \mathrm{e}^{- \mu t}
\sum_{j=0}^{k-1} \pmatrix{{k-1}\cr{j}} ( -1 )^j \mathrm{e}^{- \lambda j
t}\,\mathrm{d}t
\nonumber
\\[-8pt]
\\[-8pt]
\nonumber
&=& \int_0^\infty \mathrm{e}^{- \mu t} \{ \mathrm{e}^{- \lambda t} ( 1-\mathrm{e}^{-
\lambda t} )^{k-1} \} \,\mathrm{d}t
\end{eqnarray}
and this is an alternative derivation of the Yule--Furry linear birth
process distribution.
\begin{rem}
An alternative form of the distribution (\ref{Pk}) can be derived by
explicitly writing the Mittag--Leffler
function and conveniently manipulating the double sums obtained.
We therefore
have
%
%e3.10 ###
\begin{eqnarray}
\label{alternative}
p_k^\nu( t )
& = &\sum_{m=0}^{k-1} \frac{( - \lambda t^\nu)^m }{\Gamma
(\nu m + 1 ) } \sum_{j=0}^{k-1} \pmatrix{{k-1}\cr{j}} ( -1
)^j ( j+1
)^m\nonumber\\
&&{}+ \sum_{m=k}^\infty\frac{ ( - \lambda t^\nu)^m}{\Gamma
( \nu m +1 )}
\sum_{j=0}^{k-1} \pmatrix{{k-1}\cr{j}} ( -1 )^j ( j+1
)^m \\
& =& \frac{( \lambda t^\nu)^{k-1} ( k-1
)!}{\Gamma( \nu(
k-1 ) + 1 ) } + \sum_{m=k}^{\infty} \frac{( - \lambda
t^\nu)^m }{\Gamma
(\nu m + 1 ) } \sum_{j=0}^{k-1} \pmatrix{{k-1}\cr{j}} ( -1
)^j (
j+1 )^m. \nonumber
\end{eqnarray}
The last step of (\ref{alternative}) is justified by the following
formulas (see $0.154(6)$ and $0.154(5)$ on page~4 of Gradshteyn and Ryzhik
\cite{grad}):
%
%e3.12 ###
%e3.11 ###
\begin{eqnarray}
\sum_{k=0}^N ( -1 )^k \pmatrix{{N}\cr{k}} ( \alpha+ k
)^{n-1} &=& 0,\qquad
\mbox{valid for } N \geq n \geq1,\\\label{remarkable-thing}
\sum_{k=0}^n ( -1 )^k \pmatrix{{n}\cr{k}} ( \alpha+ k
)^n &=& ( -1 )^n n!   .
\end{eqnarray}
What is remarkable about (\ref{remarkable-thing}) is that the result is
independent of $\alpha$. This can be
ascertained as follows:
%
%e3.13 ###
\begin{equation}
\label{independence}
S_n^\alpha= \sum_{k=0}^n ( -1 )^k \pmatrix{{n}\cr{k}} \sum
_{r=0}^n \pmatrix{{n}\cr{r}} \alpha^r k^{n-r}
= \sum_{r=0}^n\pmatrix{{n}\cr{r}} \alpha^r \sum_{k=0}^n ( -1 )^k
\pmatrix{{n}\cr{k}} k^{n-r+1-1}.
\end{equation}
By formula 0.154(3) on page 4 of Gradshteyn and Ryzhik \cite{grad}, the inner sum in the
third member of~(\ref{independence}) equals zero for
$1 \leq n-r+1 \leq n$ (that is, for $1\leq r \leq n$). Therefore (see
formula~0.154(4) on page 4 of Gradshteyn and Ryzhik \cite{grad}),
%
%e3.14 ###
\begin{equation}
S_n^\alpha= \pmatrix{{n}\cr{0}} \alpha^0 \sum_{k=0}^n ( -1 )^k
\pmatrix{{n}\cr{k}} k^n = (
-1 )^n n!   .
\end{equation}
\end{rem}
We now provide a direct proof that the distribution (\ref{Pk}) sums to
unity. This is based on combinatorial
arguments and will subsequently be validated by resorting to the
representation of $N_\nu( t )$
as a composition of the Yule--Furry model with the random time $T_{2
\nu} ( t )$.
\begin{thm}
\label{th}
The distribution (\ref{Pk}) is such that
%
%e3.15 ###
\begin{equation}
\label{sums}
\sum_{k=1}^\infty p_k^\nu( t ) = \sum_{k=1}^\infty\sum
_{j=1}^k \pmatrix{{k-1}\cr{j-1}}
( -1 )^{j-1} E_{\nu,1} ( - \lambda j t^\nu) =
1 .
\end{equation}
\end{thm}
\begin{pf}
We start by evaluating the Laplace transform $L_\nu( \mu)$
of (\ref{sums}) as follows:
%
%e3.16 ###
\begin{eqnarray}
L_\nu( \mu) & =& \sum_{k=1}^\infty\sum_{j=1}^k
\pmatrix{{k-1}\cr{j-1}} ( -1 )^{j-1}
\frac{\mu^{\nu-1}}{\mu^\nu+ \lambda j}
\nonumber
\\[-8pt]
\\[-8pt]
\nonumber
&=& \frac{\mu^{\nu-1}}{\lambda} \sum_{k=1}^\infty\sum_{j=0}^{k-1}
\pmatrix{{k-1}\cr{j}}
( -1 )^j \frac{1}{{\mu^\nu}/{\lambda} + 1 + j} .
\end{eqnarray}
A crucial role is played here by the well-known formula (see Kirschenhofer \cite
{Kirschenhofer})
%e3.17 ###
\begin{equation}
\label{kirsch}
\sum_{k=0}^N \pmatrix{{N}\cr{k}} ( -1 )^k \frac{1}{x+k} = \frac{N!}{
x ( x+1 ) \cdots( x + N )} .
\end{equation}
Therefore,
%
%e3.18 ###
\begin{eqnarray}
\label{conclusion-sums}
L_\nu( \mu) & = &\frac{\mu^{\nu-1}}{\lambda} \sum
_{k=1}^\infty
\frac{( k-1 )!}{
( {\mu^\nu}/{\lambda} +1 ) ( {\mu^\nu}/{\lambda
} +2 )
\cdots( {\mu^\nu}/{\lambda} +k )}\nonumber\\
&=& \frac{\mu^{\nu-1}}{\lambda} \sum_{l=0}^\infty\frac{\Gamma(
l+1 )
\Gamma( {\mu^\nu}/{\lambda} +1 ) }{ \Gamma(
{\mu^\nu}/{
\lambda} + 1 + ( l+1 ) ) } \\
& =& \frac{\mu^{\nu-1}}{\lambda} \sum_{l=0}^\infty\mathrm{B} \biggl(
l+1, \frac{\mu^\nu}{
\lambda} +1 \biggr)
= \frac{\mu^{\nu-1}}{\lambda} \int_0^1 \sum_{l=0}^\infty x^l (
1-x )^{
{\mu^\nu}/{\lambda}}\, \mathrm{d}x \nonumber\\
& =& \frac{\mu^{\nu-1}}{\lambda} \int_0^1 ( 1-x )^{{\mu
^\nu}/{\lambda}
-1 }\, \mathrm{d}x
= \int_0^\infty \mathrm{e}^{- \mu t}\, \mathrm{d}t , \nonumber
\end{eqnarray}
where $\mathrm{B} ( p,q ) = \int_0^1 x^{p-1} ( 1-x
)^{q-1} \,\mathrm{d}x$ for $p,q > 0$.
This concludes the proof of (\ref{sums}).
\end{pf}

The presence of alternating sums in (\ref{Pk}) imposes the check that
$p_k^\nu( t )
\geq0$ for all $k$. This is the purpose of the next remark.
\begin{rem}
In order to check the non-negativity of (\ref{Pk}), we exploit the
results of the proof of Theorem \ref{th},
suitably adapted. The expression
%
%e3.19 ###
\begin{equation}
\sum_{k=1}^\infty\int_0^\infty \mathrm{e}^{- \mu t} p_k^\nu( t )
\,\mathrm{d}t = \frac{\mu^{\nu-1}}{\lambda}
\sum_{k=1}^\infty\mathrm{B} \biggl( k, \frac{\mu^\nu}{\lambda} +1 \biggr)
\end{equation}
which emerges from (\ref{conclusion-sums}) permits us to write
%
%e3.20 ###
\begin{eqnarray}
 \int_0^\infty \mathrm{e}^{- \mu t} p_k^\nu( t )\, \mathrm{d}t
&=& \int_0^1 x^{k-1} \frac{\mu^{\nu-1}}{\lambda} ( 1-x
)^{{\mu^\nu}/{\lambda}}\,\mathrm{d}x\nonumber\\
&=& \int_0^1 x^{k-1} \frac{\mu^{\nu-1}}{\lambda} \mathrm{e}^{({\mu^\nu
}/{\lambda}) \operatorname{ln} ( 1-x )}
\,\mathrm{d}x
\nonumber
\\[-8pt]
\\[-8pt]
\nonumber
& =& \int_0^1 x^{k-1} \frac{\mu^{\nu-1}}{\lambda} \mathrm{e}^{-{\mu^\nu
}/{\lambda}
 \sum _{r=1}^\infty
{x^r}/{r} }\, \mathrm{d}x\\
&=& \int_0^1 x^{k-1} \frac{ \mu^{\nu-1}}{\lambda} \mathrm{e}^{ - {\mu^\nu
x}/{\lambda} }
 \prod _{r=2}^\infty \mathrm{e}^{- { \mu^\nu x^r}/{(\lambda
r)} }\, \mathrm{d}x . \nonumber
\end{eqnarray}
The terms
%
%e3.21 ###
\begin{equation}
\mathrm{e}^{ - { \mu^\nu x^r}/{(\lambda r)} } = \mathbb{E} \mathrm{e}^{ -\mu X_r} =
\int_0^\infty \mathrm{e}^{ - \mu t}
q_\nu^r ( x, t ) \,\mathrm{d}t
\end{equation}
are the Laplace transforms of stable random variables $ X_r = S (
\sigma_r, 1, 0 ) $, where
$
\sigma_r =\break ( \frac{x^r}{\lambda r} \cos\frac{ \uppi\nu}{2}
)^{{1}/{\nu}}
$
(for details on this point, see Samorodnitsky and Taqqu \cite{taqqu}, page 15). The term
$
\frac{\mu^{\nu-1}}{2 \lambda} \exp( - \frac{ \mu^\nu |x|}{\lambda} )
$
is the Laplace transform of the solution of the fractional diffusion equation
%
%e3.22 ###
\begin{equation}
\label{diffusion}
\cases{
\displaystyle\frac{\partial^{2 \nu} u}{\partial t^{2 \nu}} = \lambda^2 \frac
{\partial^2 u}{\partial x^2},
&\quad  $0 < \nu\leq1,$ \vspace*{2pt}\cr
u ( x, 0 ) = \delta( x ), }
\end{equation}
with the additional condition that $u_t (x,0 ) = 0$ for $1/2
< \nu\leq1$,
and can be written as
%
%e3.23 ###
\begin{equation}
u_{2 \nu} ( x,t ) = \frac{1}{2 \lambda\Gamma( 1- \nu
)} \int_0^t
\frac{p_\nu( x, s )}{( t-s )^\nu} \,\mathrm{d}s
\end{equation}
(see formula (3.5) in Orsingher and Beghin \cite{ors2004}), where $p_\nu( x,1 )
= q_\nu^1 ( x,1 )$
is the stable law with
$
\sigma_1 = ( \frac{x}{\lambda} \cos\frac{\uppi\nu}{2}
)^{{1}/{\nu}} .
$
We can represent the product
%
%e3.24 ###
\begin{equation}
\label{product}
\frac{\mu^{\nu-1}}{\lambda} \mathrm{e}^{- {x \mu^\nu}/{\lambda}} \prod
_{r=2}^\infty
\mathrm{e}^{- {\mu^\nu x^r}/{(\lambda r)} } =
\int_0^\infty \mathrm{e}^{- \mu t} \biggl\{ \int_0^t u_{2 \nu} ( x,s
) q_\nu(
x, t-s )\, \mathrm{d}s \biggr\}\, \mathrm{d}t ,
\end{equation}
where
%
%e3.25 ###
\begin{equation}
\int_0^\infty \mathrm{e}^{- \mu t} q_\nu( x,t )\, \mathrm{d}t = \prod
_{r=2}^\infty \mathrm{e}^{-
{\mu^\nu x^r}/{(\lambda r)} } .
\end{equation}
Thus $q_\nu( x,t ) $ appears as an infinite convolution of
stable laws whose
parameters depend on $r$ and $x$.
In the light of (\ref{product}), we therefore have that
%
%e3.26 ###
\begin{equation}
\int_0^\infty \mathrm{e}^{- \mu t} p_k^\nu( t )\, \mathrm{d}t
= 2 \int_0^\infty \mathrm{e}^{- \mu t} \int_0^1 x^{k-1} \int_0^t u_{2 \nu}
( x,s )
q_\nu( x,t-s ) \,\mathrm{d}s \, \mathrm{d}x \, \mathrm{d}t .
\end{equation}
Since $p_k^\nu( t )$ appears as the result of the integral
of probability densities, we can conclude that
$p_k^\nu( t ) \geq0$ for all $k \geq1$ and $t > 0$.
\end{rem}

We provide an alternative proof of the non-negativity of $p_k^\nu
( t )$, $t > 0$, and
of $\sum_k p_k^\nu( t ) = 1$, based on the representation
of the fractional linear birth process
$N_\nu( t )$ as
%
%e3.27 ###
\begin{equation}
\label{iteration}
N_\nu( t ) = N ( T_{2 \nu} ( t ) ),\qquad
 0< \nu\leq1,
\end{equation}
where $T_{2 \nu} ( t )$ possesses a distribution coinciding
with the folded solution of the fractional
diffusion equation
%
%e3.28 ###
\begin{equation}
\cases{
\dfrac{\partial^{2 \nu} u}{\partial t^{2 \nu}} = \dfrac{\partial^2
u}{\partial x^2},
& \quad $0 < \nu\leq1,$\vspace*{2pt}\cr
u ( x,0 ) = \delta( x ),}
\end{equation}
with the further condition that $u_t (x,0 ) = 0$ for $1/2 <
\nu\leq1$.
\begin{thm}
The probability generating function $G_\nu( u,t ) = \mathbb{E}
u^{N_\nu( t )}$ of $N_\nu( t )$, $t>0$,
has the Laplace transform
%
%e3.29 ###
\begin{equation}
\label{lapl-trans}
\int_0^\infty \mathrm{e}^{- \mu t } G_\nu( u,t ) \,\mathrm{d}t = \int_0^\infty
\frac{u \mathrm{e}^{- \lambda t}}{
1- u ( 1- \mathrm{e}^{- \lambda t} )} \mu^{\nu-1} \mathrm{e}^{- \mu^\nu t}
\,\mathrm{d}t .
\end{equation}
\end{thm}
\begin{pf}
We evaluate the Laplace transform (\ref{lapl-trans}) as follows:
%
%e3.30 ###
\begin{eqnarray}
 \int_0^\infty \mathrm{e}^{- \mu t } G_\nu( u,t ) \,\mathrm{d}t& =&
\int_0^\infty \mathrm{e}^{- \mu t} \sum_{k=1}^\infty u^k \sum_{j=1}^k
\pmatrix{{k-1}\cr{j-1}}
( -1 )^{j-1} E_{\nu,1} ( -\lambda j t^\nu)\, \mathrm{d}t \nonumber\\
& =& \sum_{k=1}^\infty u^k \sum_{j=1}^k \pmatrix{{k-1}\cr{j-1}} ( -1
)^{j-1}
\frac{\mu^{\nu-1}}{\mu^\nu+ \lambda j}\nonumber\\
&=& \frac{\mu^{\nu-1}}{\lambda} \sum_{k=1}^\infty u^k \sum_{j=0}^{k-1}
\pmatrix{{k-1}\cr{j}}
( -1 )^j \frac{1}{{\mu^\nu}/{\lambda} + 1 + j} \qquad
  \mbox{(by (\ref{kirsch}))} \nonumber\\
& =& \frac{\mu^{\nu-1}}{\lambda} \sum_{k=1}^\infty u^k \frac{ (
k-1 )!}{
( {\mu^\nu}/{\lambda} +1 ) ( {\mu^\nu}/{\lambda
} +2 )
\cdots( {\mu^\nu}/{\lambda} + k ) }\nonumber\\
&=& \frac{u \mu^{\nu-1}}{\lambda} \sum_{l=0}^\infty u^l \frac{l!}{
( {\mu^\nu}/{\lambda}
+1 ) \cdots({\mu^\nu}/{\lambda} +1+l ) }\\
& =& \frac{ u \mu^{\nu-1}}{\lambda} \sum_{l=0}^\infty u^l \mathrm{B}
\biggl( l+1, \frac{\mu^\nu}{\lambda}
+1 \biggr)\nonumber\\
&=& \frac{ u \mu^{\nu-1}}{\lambda} \int_0^1 \sum_{l=0}^\infty u^l x^l
( 1-x
)^{{\mu^\nu}/{\lambda}} \,\mathrm{d}x \qquad \mbox{(for $ 0 < ux < 1$)}
\nonumber\\
& =& \frac{ u \mu^{\nu-1}}{\lambda} \int_0^1 \frac{ ( 1-x
)^{{\mu^\nu}/{\lambda}
}}{ ( 1 - ux ) } \,\mathrm{d}x =   ( 1-x = \mathrm{e}^{-\lambda t} )\nonumber\\
&=& \int_0^\infty\frac{u \mathrm{e}^{- \lambda t}}{1- u ( 1- \mathrm{e}^{- \lambda
t} )} \mathrm{e}^{-t \mu^\nu}
\mu^{\nu-1}\, \mathrm{d}t \nonumber.
\end{eqnarray}
\upqed\end{pf}
\begin{rem}
In order to extract from (\ref{lapl-trans}) the representation (\ref
{iteration}), we note that
%
%e3.31 ###
\begin{eqnarray}
&& \int_0^\infty \mathrm{e}^{- \mu t} \Biggl\{ \sum_{k=0}^\infty u^k \operatorname{Pr}
\{ N ( T_{2 \nu}
( t ) ) = k \} \Biggr\} \,\mathrm{d}t \nonumber\\
&&\quad = \int_0^\infty \mathrm{e}^{ - \mu t} \Biggl\{ \int_0^\infty\sum_{k=0}^\infty
u^k \operatorname{Pr} \{
N ( s ) = k \} f_{T_{2 \nu}} ( s,t )\, \mathrm{d}s
\Biggr\}\, \mathrm{d}t \\
&&\quad=
\int_0^\infty G ( u,s ) \mu^{\nu-1} \mathrm{e}^{- \mu^\nu s}\, \mathrm{d}s,
\nonumber
\end{eqnarray}
which coincides with (\ref{lapl-trans}). It can be shown that
%
%e3.32 ###
\begin{equation}
\label{lapl-trans-folded}
\int_0^\infty \mathrm{e}^{- \mu t } f_{T_{2 \nu}} ( s,t )\, \mathrm{d}t = \mu
^{\nu-1} \mathrm{e}^{- s \mu^\nu},\qquad
  s >0,
\end{equation}
is the Laplace transform of the folded solution to
%
%e3.33 ###
\begin{equation}
\label{folding}
\dfrac{\partial^{2 \nu} u}{\partial t^{2 \nu}} = \dfrac{\partial^2
u}{\partial s^2},\qquad
0 < \nu\leq1,
\end{equation}
with the initial condition $u ( s,0 ) = \delta( s
)$ for
$0 < \nu\leq1$ and also $u_t ( s,0 ) =0$ for $1/2 < \nu
\leq1$.

In the light of (\ref{iteration}), the non-negativity of $p_k^\nu
( t )$ is immediate because
%
%e3.34 ###
\begin{equation}
\label{iter-relat}
\operatorname{Pr} \{ N_\nu( t ) = k \} = \int_0^\infty
\operatorname{Pr} \{ N
( s ) = k \} \operatorname{Pr} \{ T_{2 \nu} ( t
) \in \mathrm{d}s \} .
\end{equation}
The relation (\ref{iter-relat}) immediately leads to the conclusion that
$
\sum_{k=1}^\infty\operatorname{Pr} \{ N_\nu( t ) = k
\} = 1 .
$
\end{rem}

Some explicit expressions for (\ref{iter-relat}) can be given when the
$\operatorname{Pr} \{ T_{2 \nu} ( t )
\in \mathrm{d}s \}$ can be worked out in detail.

We know that for $\nu= 1/2^n$, we have that
%
%e3.35 ###
\begin{eqnarray}
\label{iterated-distr2}
&&\operatorname{Pr} \{ T_{{1}/{2^{n-1}}} ( t ) \in \mathrm{d}s
\}\nonumber\\
&&\quad=\operatorname{Pr} \{ \vert\mathcal{B}_1 ( \vert\mathcal{B}_2
( \cdots\vert\mathcal{B}_n
( t ) \vert\cdots) \vert )
\in \mathrm{d}s \}\\
&&\quad =\mathrm{d}s 2^n \int_0^\infty\frac{\mathrm{e}^{- {s^2}/{(4 \omega_1)}}}{\sqrt{4 \uppi
\omega_1}} \,\mathrm{d} \omega_1
\int_0^\infty\frac{\mathrm{e}^{- {\omega_1^2}/{(4 \omega_2)}}}{\sqrt{4 \uppi
\omega_2}}\, \mathrm{d} \omega_2 \cdots
\int_0^\infty\frac{\mathrm{e}^{- {\omega_{n-1}^2}/{(4t)}}}{\sqrt{4 \uppi t}} \,\mathrm{d}
\omega_{n-1} . \nonumber
\end{eqnarray}
For details concerning (\ref{iterated-distr2}), see Theorem $2.2$ of Orsingher and Beghin
\cite{ors2008},
where the differences of the constants depend on the fact that the
diffusion coefficient in equation (\ref{folding})
equals $1$ instead of $2^{( 1/ 2^n ) -2}$.
The distribution (\ref{iterated-distr2}) represents the density of the
folded $( n-1 )$-times iterated
Brownian motion and therefore $\mathcal{B}_1, \ldots, \mathcal{B}_n$
are independent Brownian motions with
volatility equal to $2$.

For $\nu= 1/3$, the process (\ref{iteration}) has the form
$
N_{{1}/{3}} ( t ) = N ( \vert\mathcal{A} (
t ) \vert)
$,
where $\mathcal{A} ( t )$ is a process whose law is the
solution of
%
%e3.36 ###
\begin{equation}
\label{diff-eq-airy}
\frac{\partial^{{2}/{3}} u}{\partial t^{{2}/{3}}} = \frac
{\partial^2 u}{\partial x^2},\qquad
  u ( x,0 ) = \delta( x ) .
\end{equation}
In Orsingher and Beghin \cite{ors2008}, it is shown that the solution to (\ref
{diff-eq-airy}) is
%
%e3.37 ###
\begin{equation}
u_{{2}/{3}} ( x,t ) = \frac{3}{2} \frac{1}{\sqrt[3]{3
t}} \mathcal{A}_i \biggl(
\frac{\vert x \vert}{\sqrt[3]{3t}} \biggr),
\end{equation}
where
%
%e3.38 ###
\begin{equation}
\mathcal{A}_i ( x ) = \frac{1}{\uppi} \int_0^\infty\cos
\biggl( \alpha x + \frac{\alpha^3}{3}
\biggr)\, \mathrm{d} \alpha
\end{equation}
is the Airy function. Therefore, in this case, the distribution (\ref
{iter-relat}) has the form
%
%e3.39 ###
\begin{equation}
p_k^{{1}/{3}} ( t ) = \int_0^\infty \mathrm{e}^{- \lambda s}
( 1- \mathrm{e}^{- \lambda s} )^{k-1}
\frac{3}{\sqrt[3]{3t}} \mathcal{A}_i \biggl( \frac{ s }{
\sqrt[3]{3t}} \biggr)\, \mathrm{d}s,\qquad  k \geq1, t > 0.
\end{equation}
\begin{rem}
\label{rm}
From (\ref{diff-diff-eq}), it is straightforward to show that the
probability generating function
$G_\nu( u,t ) = \mathbb{E} u^{N_\nu( t )}$
satisfies the partial differential
equation
%
%e3.40 ###
\begin{equation}
\label{frac-part-g}
\cases{
\dfrac{\partial^\nu}{\partial t^\nu} G ( u,t ) = \lambda u
( u-1 )
\dfrac{\partial}{\partial u} G ( u,t ), &\quad  $0 < \nu\leq1,$\vspace*{2pt}\cr
G ( u,0 ) = u, &}
\end{equation}
and thus
$
\mathbb{E} N_\nu( t ) =  \frac{\partial G}{\partial
u} |_{u=1}
$
is the solution to
%
%e3.41 ###
\begin{equation}
\label{eq-expect}
\cases{
\dfrac{ \mathrm{d}^\nu}{\mathrm{d} t^\nu} \mathbb{E} N_\nu= \lambda
\mathbb{E} N_\nu,&
\quad $0 < \nu\leq1,$ \vspace*{2pt}\cr
\mathbb{E} N_\nu( 0 ) = 1. &}
\end{equation}
The solution of (\ref{eq-expect}) is
%
%e3.42 ###
\begin{equation}
\label{result-expect}
\mathbb{E} N_\nu( t ) = E_{\nu,1} ( \lambda t^\nu
), \qquad  t>0 .
\end{equation}
Clearly, the result (\ref{result-expect}) can be also derived by
evaluating the
Laplace transform
\begin{eqnarray*}
 \int_0^\infty \mathrm{e}^{- \mu t} \mathbb{E} N_\nu( t ) \,\mathrm{d}t
 &=&
\int_0^\infty \mathrm{e}^{- \mu t} \Biggl\{
\sum_{k=1}^\infty k \int_0^\infty\operatorname{Pr} \{ N ( s )
= k \}
\operatorname{Pr} \{ T_{2 \nu} ( t ) \in \mathrm{d}s \}
\Biggr\} \,\mathrm{d}t \\
& =& \int_0^\infty \mathrm{e}^{- \mu t} \int_0^\infty \mathrm{e}^{\lambda s} \operatorname{Pr}
\{
T_{2 \nu} ( t ) \in \mathrm{d}s \}\, \mathrm{d}t\\
&= &\int_0^\infty \mathrm{e}^{ \lambda s} \mu^{\nu-1} \mathrm{e}^{- s \mu^\nu} \mathrm{d}s
= \frac{\mu^{\nu-1}}{\mu^\nu- \lambda} = \int_0^\infty \mathrm{e}^{ - \mu t}
E_{\nu,1} ( \lambda
t^\nu) \,\mathrm{d}t
\end{eqnarray*}

and this verifies (\ref{result-expect}). The mean value (\ref
{result-expect}) can be obtained in a third
manner:
%
%e3.43 ###
\begin{eqnarray}
 \int_0^\infty \mathrm{e}^{- \mu t} \mathbb{E} N_\nu( t ) &=& \sum
_{k=1}^\infty k \sum_{j=1}^k
\pmatrix{{k-1}\cr{j-1}} ( -1 )^{j-1} \int_0^\infty E_{\nu,1}
( - \lambda j t^\nu)
\mathrm{e}^{- \lambda t}\, \mathrm{d}t \nonumber\\
& =& \sum_{k=1}^\infty k \sum_{j=1}^k \pmatrix{{k-1}\cr{j-1}} ( -1
)^{j-1} \frac{\mu^{\nu-1}}{
\mu^\nu+ \lambda j}\nonumber\\
&= &\frac{\mu^{\nu-1}}{\lambda} \sum_{k=1}^\infty k \sum_{j=0}^{k-1}
\pmatrix{{k-1}\cr{j}} ( -1
)^j \frac{1}{{\mu^\nu}/{\lambda} + 1 + j}
\nonumber
\\[-8pt]
\\[-8pt]
\nonumber
& = &\frac{\mu^{\nu-1}}{\lambda} \sum_{k=1}^\infty k \frac{ ( k-1
)!}{({\mu^\nu}/{\lambda} +1 ) \cdots( {\mu^\nu}/{\lambda
} + k ) } \\
&=&
\frac{\mu^{\nu-1}}{\lambda} \sum_{k-1}^\infty k \frac{ \Gamma(
k ) \Gamma
( {\mu^\nu}/{\lambda} + 1 )}{ \Gamma({\mu^\nu}/{\lambda}
+ k +1 ) }\nonumber\\
& =& \frac{\mu^{\nu-1}}{\lambda} \int_0^1 \sum_{k=1}^\infty k x^{k-1}
( 1-x )^{{\mu^\nu}/{
\lambda}}
= \frac{\mu^{\nu-1}}{
\mu^\nu- \lambda}
= \int_0^\infty \mathrm{e}^{- \mu t} E_{\nu,1} ( \lambda t^\nu)\, \mathrm{d}t.
\nonumber
\end{eqnarray}
\end{rem}
The result of Remark \ref{rm}, $\mathbb{E} N_\nu( t ) =
E_{\nu,1} (
\lambda t^\nu)$, should be compared with the results of Uchaikin, Cahoy and Sibatov \cite{sib}.
%
%f1 ###
\begin{figure}[b]

\includegraphics{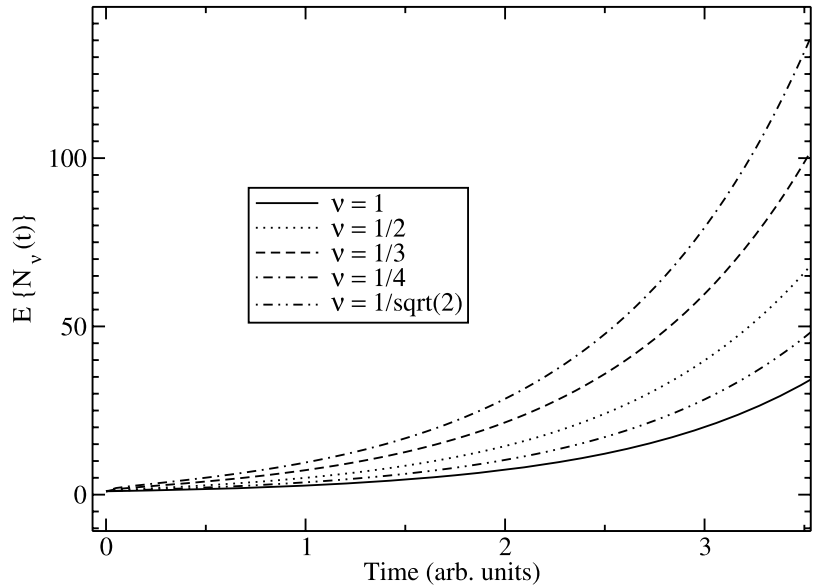}

\caption{Mean number of individuals at time $t$ for
various values of $\nu$.}\label{expectation}
\end{figure}

An interesting representation of (\ref{result-expect}) following from
(\ref{iteration}) gives that
%
%e3.44 ###
\begin{equation}
\mathbb{E} N_\nu( t ) = \int_0^\infty \mathrm{e}^{\lambda s} \operatorname
{Pr} \{ T_{2 \nu} ( t )
\in \mathrm{d}s \} = \int_0^\infty\mathbb{E} N ( s ) \operatorname
{Pr} \{ T_{2 \nu} ( t )
\in \mathrm{d}s \} .
\end{equation}
The expansion of the population subject to the law of the fractional
birth process is increasingly rapid as the
order of fractionality $\nu$ decreases. This is shown in Figure \ref
{expectation} and this behavior is due to
the increasing structure of the gamma function\ for $\nu>0$ appearing
in the Mittag--Leffler function~$E_{\nu,1}$.
This qualitative feature of the process being investigated here shows
that it conveniently applies to
explosively expanding populations.
\begin{rem}
By twice deriving (\ref{frac-part-g}) with respect to $u$, we obtain
the fractional equation for the second-order
factorial moment
%
%e3.45 ###
\begin{equation}
\mathbb{E} \bigl\{ N_\nu( t ) \bigl( N_\nu( t
) - 1 \bigr) \bigr\}
= g_\nu( t ),
\end{equation}
that is,
%
%e3.46 ###
\begin{equation}
\label{factorial}
\cases{
\dfrac{\partial^\nu}{\partial t^\nu} g_\nu( t ) = 2
\lambda g_\nu( t )
+ 2 \lambda\mathbb{E} N_\nu( t ), &\quad  $0<\nu\leq1,$\vspace*{2pt}\cr
g_\nu( 0 ) = 0.&}
\end{equation}
The Laplace transform of the solution to (\ref{factorial}) is
%
%e3.47 ###
\begin{eqnarray}
\label{laplace-factorial}
H_\nu( t ) &= &\int_0^\infty \mathrm{e}^{- \mu t} g_\nu( t
)\, \mathrm{d}t =
\frac{2 \lambda\mu^{\nu-1} }{ ( \mu^\nu- \lambda)
( \mu^\nu- 2 \lambda)}
\nonumber
\\[-8pt]
\\[-8pt]
\nonumber
&= &2 \mu^{\nu-1} \biggl\{ \frac{1}{ \mu^\nu- 2 \lambda} - \frac{1}{ \mu
^\nu- \lambda} \biggr\} .
\end{eqnarray}
The inverse Laplace transform of (\ref{laplace-factorial}) is
%
%e3.48 ###
\begin{equation}
\label{inverse-lapl-factorial}
\mathbb{E} \bigl\{ N_\nu( t ) \bigl( N_\nu( t
) -1 \bigr) \bigr\} =
2 E_{\nu,1} ( 2 \lambda t^\nu) - 2 E_{\nu,1} (
\lambda t^\nu) .
\end{equation}
It is now straightforward to obtain the variance from (\ref
{inverse-lapl-factorial}),
%e3.49 ###
\begin{equation}
\label{variance}
\operatorname{\mathbb{V}ar} N_\nu( t ) = 2 E_{\nu,1} ( 2
\lambda t^\nu)
- E_{\nu,1} ( \lambda t^\nu) - E_{\nu,1}^2 ( \lambda
t^\nu) .
\end{equation}
For $\nu=1$, we retrieve from (\ref{variance}) the well-known
expression of the variance of the linear
birth process
%
%e3.50 ###
\begin{equation}
\operatorname{\mathbb{V}ar} N_1 ( t ) = \mathrm{e}^{\lambda t } (
\mathrm{e}^{\lambda t} -1 ) .
\end{equation}
\end{rem}
\begin{rem}
If $X_1, \ldots, X_n$ are i.i.d.~random variables with common
distribution $F ( x ) =\break \operatorname{Pr} ( X < x )$,
then we can write the following probability:
%
%e3.51 ###
\begin{eqnarray}
&& \operatorname{Pr} \bigl\{ \max\bigl( X_1, \ldots, X_{N_\nu( t
)} \bigr) < x \bigr\}\nonumber\\
&&\quad =
\sum_{k=1}^\infty( \operatorname{Pr} \{ X < x \} )^k
\operatorname{Pr} \{ N_\nu( t ) =
k \} \qquad   \mbox{(by (\ref{iteration}))}
\nonumber
\\[-8pt]
\\[-8pt]
\nonumber
&&\quad = \int_0^\infty G ( F ( x ) , s ) \operatorname{Pr}
\{ T_{2 \nu} ( t ) \in \mathrm{d}s
\}\\
&&\quad =
\int_0^\infty\frac{F ( x ) \mathrm{e}^{- \lambda s}}{1 - F ( x
) ( 1 - \mathrm{e}^{- \lambda s}
) } \operatorname{Pr} \{ T_{2 \nu} ( t ) \in \,\mathrm{d}s \}
. \nonumber
\end{eqnarray}
Analogously, we have that
%
%e3.52 ###
\begin{eqnarray}
&& \operatorname{Pr} \bigl\{ \min\bigl( X_1, \ldots, X_{N_\nu( t
)} \bigr) > x \bigr\}
\nonumber
\\[-8pt]
\\[-8pt]
\nonumber
&&\quad= \int_0^\infty\frac{ ( 1 - F ( x ) ) \mathrm{e}^{-
\lambda s} }{1 - ( 1 -F ( x )
) ( 1 - \mathrm{e}^{- \lambda s} ) } \operatorname{Pr} \{ T_{2 \nu
} ( t ) \in \,\mathrm{d}s \} .
\end{eqnarray}
\end{rem}

\begin{rem}
If the initial number of components of the population is $n_0$, then
the p.g.f.~becomes
%
%e3.53 ###
\begin{eqnarray}
\label{pgf}
&&\mathbb{E} \bigl( u^{N_\nu( t )} | N_\nu
( 0 ) = n_0 \bigr)
\nonumber
\\[-8pt]
\\[-8pt]
\nonumber
&&\quad=\sum_{k=0}^\infty u^{k+n_0} \int_0^\infty \mathrm{e}^{- \lambda z n_0}
\pmatrix{{n_0 + k -1}\cr{k}} (
1 - \mathrm{e}^{- \lambda z} )^k \operatorname{Pr} \{ T_{2 \nu} ( t
) \in \,\mathrm{d}z \} .
\end{eqnarray}
From (\ref{pgf}), we can extract the distribution of the population
size at time $t$ as
%
%e3.54 ###
\begin{eqnarray}
\label{distr-size}
&& \operatorname{Pr} \{ N_\nu( t ) = k+ n_0 |
N_\nu( 0 ) = n_0 \}
\nonumber
\\[-8pt]
\\[-8pt]
\nonumber
&&\quad = \pmatrix{{n_0 + k -1}\cr{k}} \int_0^\infty \mathrm{e}^{- \lambda z n_0} ( 1 -
\mathrm{e}^{- \lambda z} )^k
\operatorname{Pr} \{ T_{2 \nu} ( t ) \in\, \mathrm{d}z \},\qquad    k
\geq0 .
\end{eqnarray}
If we write $k + n_0 =k'$, then we can rewrite (\ref{distr-size}) as
%
%e3.55 ###
\begin{eqnarray}
\label{distr-size-other}
&& \operatorname{Pr} \{ N_\nu( t ) = k' | N_\nu
( 0 ) = n_0 \}
\nonumber
\\[-8pt]
\\[-8pt]
\nonumber
&&\quad= \pmatrix{{k' -1}\cr{k' - n_0}} \int_0^\infty \mathrm{e}^{- \lambda z n_0} ( 1
- \mathrm{e}^{- \lambda z} )^{k'-n_0}
\operatorname{Pr} \{ T_{2 \nu} ( t ) \in \,\mathrm{d}z \},\qquad
k' \geq n_0,
\end{eqnarray}
where $k'$ is the number of individuals in the population at time $t$.
For $n_0 = 1$, formulas (\ref{distr-size}),
(\ref{distr-size-other}) coincide with (\ref{Pk}).
The random time $T_{2 \nu} ( t )$, $t>0$, appearing in (\ref
{distr-size}) and (\ref{distr-size-other})
has a distribution which is related to the fractional equation
%
%e3.56 ###
\begin{equation}
\frac{\partial^{2 \nu} u}{\partial t^{2 \nu}} = \frac{\partial^2
u}{\partial z^2},\qquad
0 < \nu\leq1 .
\end{equation}
It is possible to slightly change the structure of formulas (\ref
{distr-size}) and (\ref{distr-size-other})
by means of the transformation $\lambda z = y$ so that the
distribution of $T_{2 \nu} ( t )$ becomes
related to the equation
%
%e3.57 ###
\begin{equation}
\label{equa}
\frac{\partial^{2 \nu} u}{\partial t^{2 \nu}} = \lambda^2 \frac
{\partial^2 u}{\partial y^2},\qquad
0 < \nu\leq1,
\end{equation}
where \eqref{linear-rates} shows the connection between the diffusion
coefficient in \eqref{equa} and the birth rate.
\end{rem}
\begin{rem}
If we assume that the initial number of individuals in the population
is $N_\nu( 0 ) = n_0$, then
we can generalize the result (\ref{Pk}) offering a representation of
the distribution of $N_\nu( t )$
alternative to (\ref{distr-size-other}). If we take the Laplace
transform of (\ref{distr-size-other}), then we have
that
\begin{eqnarray}
\label{gener-pk}
&& \int_0^\infty \mathrm{e}^{- \mu t} \operatorname{Pr} \{ N_\nu( t ) =
k + n_0  | N_\nu
( 0 ) = n_0 \} \,\mathrm{d}t \nonumber\\
&&\quad = \int_0^\infty\pmatrix{{n_0 + k -1}\cr{k}} \int_0^\infty \mathrm{e}^{- \lambda z
n_0} ( 1- \mathrm{e}^{- \lambda z}
)^k \operatorname{Pr} \{ T_{2 \nu} ( t ) \in\, \mathrm{d}z \}\,
\mathrm{d}t \qquad   \operatorname{(by (\ref{lapl-trans-folded}))} \nonumber\\
&&\quad = \int_0^\infty\pmatrix{{n_0 + k -1}\cr{k}} \mathrm{e}^{- \lambda z n_0} ( 1 -
\mathrm{e}^{- \lambda z} )^k
\mu^{\nu-1} \mathrm{e}^{- \mu^\nu z}\, \mathrm{d}z
\nonumber
\\[-8pt]
\\[-8pt]
\nonumber
&&\quad= \pmatrix{{n_0 + k -1}\cr{k}} \mu^{\nu-1} \int_0^\infty \mathrm{e}^{- z (
\lambda n_0 + \mu^\nu)}
( 1 - \mathrm{e}^{- \lambda z} )^k\, \mathrm{d}z \nonumber\\
&&\quad= \pmatrix{{n_0 + k -1}\cr{k}} \mu^{\nu-1} \sum_{r=0}^k
\pmatrix{{k}\cr{r}}
( -1 )^r \int_0^\infty
\mathrm{e}^{- z ( \lambda n_0 + \lambda r + \mu^\nu)} \,\mathrm{d}z \nonumber
\nonumber\\
&&\quad = \pmatrix{{n_0 + k -1}\cr{k}} \mu^{\nu-1} \sum_{r=0}^k \pmatrix{{k}\cr{r}} (
-1 )^r \frac{1}{ \lambda
( n_0 + r ) + \mu^\nu}  . \nonumber
\end{eqnarray}
By taking the inverse Laplace transform of (\ref{gener-pk}), we have that
%
%e3.58 ###
\begin{eqnarray}
\label{gener-pk2}
&&\operatorname{Pr} \{ N_\nu( t ) = k + n_0  |
N_\nu( 0 ) = n_0
\}
\nonumber
\\[-8pt]
\\[-8pt]
\nonumber
&&\quad= \pmatrix{{n_0 + k -1}\cr{k} }\sum_{r=0}^k \pmatrix{{k}\cr{r}} ( -1 )^r
E_{\nu,1} \bigl(
- ( n_0 + r ) \lambda t^\nu\bigr).
\end{eqnarray}
From (\ref{gener-pk2}), we can infer the interesting information
%
%e3.59 ###
\begin{eqnarray}
&& \operatorname{Pr} \{ N_\nu( \mathrm{d}t ) = n_0 +1 |
N_\nu( 0 ) = n_0 \}\nonumber\\
&&\quad=
n_0 \sum_{r=0}^1 \pmatrix{{1}\cr{r}} ( -1 )^r E_{\nu,1} \bigl( -
( n_0 + r )
\lambda( \mathrm{d}t )^\nu\bigr) \\
&&\quad = n_0 \bigl[ E_{\nu,1} \bigl( - n_0 \lambda( \mathrm{d}t ) ^\nu
\bigr) - E_{\nu,1}
\bigl( - \lambda( n_0 + 1 ) ( \mathrm{d}t )^\nu\bigr)
\bigr]
\sim n_0 \frac{\lambda( \mathrm{d}t )^\nu}{\Gamma( \nu+ 1
)} \nonumber
\end{eqnarray}
by writing only the lower order terms.
This shows that the probability of a new offspring at the beginning of
the process is proportional
to $( \mathrm{d}t )^\nu$ and to the initial number of progenitors.
From our point of view, this is the most important qualitative
feature of our results since it makes explicit the dependence on the
order $\nu$ of the fractional birth process.
\end{rem}
\begin{thm}
\label{theo}
The Laplace transform of the probability generating function $G_\nu
( t, u )$ of the fractional
linear birth process has the form
%
%e3.60 ###
\begin{equation}
\label{form-of-h}
\hspace*{-14pt}H_\nu( \mu,u ) = \int_0^\infty \mathrm{e}^{- \mu t} G_\nu(
t,u ) \,\mathrm{d}t =
\frac{ u \mu^{\nu-1}}{\lambda} \int_0^1 \frac{( 1-x
)^{{\mu^\nu}/{\lambda}}}{
1 -x u} \,\mathrm{d}x,\qquad  \mbox{$0<u<1$, $\mu>0$}.
\end{equation}
\end{thm}
\begin{pf}
We saw above that the function $G_\nu$ solves the Cauchy problem
%
%e3.61 ###
\begin{equation}
\label{cauchy-for-g}
\cases{
\dfrac{\partial^\nu G_\nu}{\partial t^\nu} = \lambda u ( u-1
)
\dfrac{\partial G_\nu}{\partial u}, &\quad  $ 0 < \nu\leq1,$\vspace*{2pt}\cr
G_\nu( u, 0 ) = u .& }
\end{equation}
By taking the Laplace transform of (\ref{cauchy-for-g}), we have that
%
%e3.62 ###
\begin{equation}
\label{lapl-for-g}
\mu^\nu H_\nu- \mu^{\nu-1} u = \lambda u ( u-1 ) \frac{
\partial H_\nu}{
\partial u} .
\end{equation}
By inserting (\ref{form-of-h}) into (\ref{lapl-for-g}) and performing
some integrations by parts, we have that
\begin{eqnarray}
&&\frac{u \mu^{2 \nu-1}}{\lambda}  \int_0^1 \frac{( 1-x
)^{{\mu^\nu}/{\lambda}}}{
1- xu} \,\mathrm{d}x - u \mu^{\nu-1}\nonumber\\
&&\quad=  \lambda u ( u-1 ) \biggl[ \frac{\mu^{\nu-1}}{\lambda
} \int_0^1 \frac{( 1-x
)^{{\mu^\nu}/{\lambda}}}{ 1- xu}\, \mathrm{d}x + \frac{u \mu^{\nu
-1}}{\lambda} \int_0^1
\frac{ ( 1-x )^{{\mu^\nu}/{\lambda}} x}{( 1-xu
)^2} \,\mathrm{d}x \biggr] \nonumber\\
&&\quad= \lambda u ( u-1 ) \biggl[ \frac{\mu^{\nu-1}}{\lambda
} \int_0^1 \frac{( 1-x)^{{\mu^\nu}/{\lambda}}}{ 1- xu} \,\mathrm{d}x + \frac{\mu^{\nu
-1}}{\lambda} \frac{x (
1-x )^{{\mu^\nu}/{\lambda}}}{1-xu} \bigg\vert_{x=0}^{x=1}
\nonumber
\\[-8pt]
\\[-8pt]
\nonumber
&&\quad \phantom{= \lambda u ( u-1 ) \biggl[}{} - \frac{\mu^{\nu-1}}{\lambda} \int_0^1 \frac{( 1-x
)^{{\mu^\nu}/{\lambda}}
}{ ( 1-xu ) }\, \mathrm{d}x + \frac{\mu^{2 \nu-1}}{\lambda^2} \int
_0^1 \frac{ x (
1 -x )^{{\mu^\nu}/{\lambda} -1}}{( 1-xu ) } \,\mathrm{d}x
\biggr]\qquad  \nonumber\\
&&\quad=  \frac{u ( u -1 ) \mu^{2 \nu-1}}{\lambda} \int_0^1
\frac{x ( 1-x
)^{{\mu^\nu}/{\lambda} -1}}{( 1 - xu )}\, \mathrm{d}x\nonumber\\
&&\quad= - u \mu^{\nu-1} + \frac{ u \mu^{2 \nu-1}}{\lambda} \int_0^1 \frac{
( 1-x )^{
{\mu^\nu}/{\lambda}}}{( 1 -xu )} \,\mathrm{d}x , \nonumber
\end{eqnarray}
and this concludes the proof of Theorem \ref{theo}.
\end{pf}
\begin{rem}
We note that $H_\nu( \mu,u ) \vert_{u=1} = 1/\mu$
because $G_\nu( t,1 ) = 1$.
Furthermore,
%
%e3.63 ###
\begin{equation}
\frac{\partial H_\nu( \mu,u )}{\partial u} \bigg\vert
_{u=1} = \frac{\mu^{\nu-1}}{
\mu^\nu- \lambda} = \int_0^\infty \mathrm{e}^{- \mu t} E_{\nu,1} (
\lambda t^\nu)\, \mathrm{d}t ,
\end{equation}
which accords well with (\ref{result-expect}).
\end{rem}

\section*{Acknowledgement}
The authors are pleased to acknowledge
the remarks of an unknown
referee which improved the quality of this paper.

% \bibliographystyle{myabbrv}
% \bibliographystyle{abbrvnat}
% \bibliographystyle{chicago}
% \bibliography{frac}

\begin{thebibliography}{19}

%b1 ###
\bibitem{bartlett}
Bartlett, M.S. (1978).
 \textit{An Introduction to Stochastic Processes, with Special
Reference to Methods and Applications}, 3rd ed.
Cambridge: Cambridge Univ. Press.
\MR{0475536}

%b2 ###
\bibitem{orsbeg}
Beghin, L. and Orsingher, E. (2009).
Fractional Poisson processes and related planar random motions.
 \textit{Electron. J. Probab.} \textbf{14} 1970--1827.
\MR{2535014}

%b3 ###
\bibitem{cahoy}
Cahoy, D.O. (2007).
Fractional Poisson processes in terms of alpha-stable
densities. Ph.D. thesis.

%b4 ###
\bibitem{feller1}
Feller, W. (1968).
 \textit{An Introduction to Probability Theory and Its Applications,
Volume 1}, 3rd ed. New York: Wiley.
\MR{0228020}


%b5 ###
\bibitem{Skoro}
Gikhman, I.I. and Skorokhod, A.V. (1996).
\textit{Introduction to the Theory of Random Processes}.
New York: Dover Publications.
\MR{1435501}

%b6 ###
\bibitem{grad}
Gradshteyn, I.S. and Ryzhik, I.M. (1980).
\textit{Table of Integrals, Series, and Products}.
New York: Academic Press.
\MR{0582453}

%b7 ###
\bibitem{jumarie}
Jumarie, G. (2001).
Fractional master equation: Non-standard analysis and
Liouville--Riemann derivative.
 \textit{Chaos Solitons Fractals} \textbf{12} 2577--2587.
\MR{1851079}

%b8 ###
\bibitem{Kirschenhofer}
Kirschenhofer, P. (1996).
A note on alternating sums.
 \textit{Electron. J. Combin.} \textbf{3} 1--10.
\MR{1392492}

%b9 ###
\bibitem{laskin}
Laskin, N. (2003). Fractional Poisson process.
\textit{Commun. Nonlinear Sci. Numer. Simul.} \textbf{8} 201--213.
\MR{2007003}

%b10 ###
\bibitem{ors2004}
Orsingher, E. and Beghin, L. (2004).
Time-fractional telegraph equations and telegraph processes with
Brownian time.
\textit{Probab. Theory Related Fields} \textbf{128} 141--160.
\MR{2027298}

%b11 ###
\bibitem{ors2008}
Orsingher, E. and Beghin, L. (2009).
Fractional diffusion equations and processes with randomly-varying
time.
\textit{Ann. Probab.} \textbf{37} 206--249.
\MR{2489164}

%b12 ###
\bibitem{podlubny}
Podlubny, I. (1999).
 \textit{Fractional Differential Equations}.
San Diego: Academic Press.
\MR{1658022}

%b13 ###
\bibitem{repin}
Repin, O.N. and Saichev, A.I. (2000).
 {Fractional Poisson law}.
 \textit{Radiophys. and Quantum Electronics} \textbf{43} 738--741.
\MR{1910034}

%b14 ###
\bibitem{taqqu}
Samorodnitsky, G. and Taqqu, M.S. (1994).
\textit{Stable Non-Gaussian Random Processes: Stochastic
Models with
Infinite Variance}.
New York: Chapman and Hall.
\MR{1280932}

%b15 ###
\bibitem{sibatov}
Uchaikin, V.V. and Sibatov, R.T. (2008).
A fractional Poisson process in a model of dispersive charge
transport in semiconductors.
\textit{Russian J. Numer. Anal. Math. Modelling} \textbf{23} 283--297.
\MR{2414873}

%b16 ###
\bibitem{sib}
Uchaikin, V.V., Cahoy, D.O. and Sibatov, R.T. (2008).
Fractional processes: From Poisson to branching one.
\textit{Int. J. Bifurcation Chaos}
\textbf{18} 2717--2725.
\MR{2479327}

%b17 ###
\bibitem{xiao1}
Wang, X.-T. and Wen, Z.-X. (2003).
Poisson fractional processes.
\textit{Chaos Solitons Fractals} \textbf{18} 169--177.
\MR{1984556}

%b18 ###
\bibitem{xiao2}
Wang, X.-T., Wen, Z.-X. and Zhang, S.-Y. (2006).
Fractional Poisson process (II).
\textit{Chaos Solitons Fractals} \textbf{28} 143--147.
\MR{2174587}

%b19 ###
\bibitem{wang}
Wang, X.-T., Zhang, S.-Y. and Fan, S. (2007).
Nonhomogeneous fractional Poisson processes.
\textit{Chaos Solitons Fractals} \textbf{31} 236--241.
\MR{2263284}

\end{thebibliography}

%

\printhistory

\end{document}